\documentclass[
leqno,          
11pt,           
a4paper         
]{amsart}

\usepackage[colorlinks,citecolor=blue,urlcolor=blue]{hyperref}          
\usepackage{courier}                                                    
\usepackage{amssymb, amsmath, amsfonts, amsthm}                         
\usepackage{mathtools, cite, enumerate, rotating, framed, lipsum}       
\usepackage{fancybox, bbm}                                              
\usepackage[dvipsnames]{xcolor}                                         
\usepackage[polish, english]{babel}                                     
\usepackage[T1]{fontenc}                                                
\usepackage[utf8]{inputenc}                                           
\usepackage[normalem]{ulem}
\usepackage{yfonts}
\usepackage{mathrsfs}
\usepackage{dsfont}
\usepackage{tikz}
\usepackage{xstring}

\DeclareMathAlphabet{\mathpzc}{OT1}{pzc}{m}{it}
 \numberwithin{equation}{section}                        

\newcommand{\thmcount}{equation}                 
\newcounter{specialcounter}

\newtheorem{Thm}[\thmcount]{Theorem}
\newtheorem{Sthm}[specialcounter]{Theorem}
\newtheorem{Cor}[\thmcount]{Corollary}
\newtheorem{Lem}[\thmcount]{Lemma}
\newtheorem{Prop}[\thmcount]{Proposition}
\newtheorem{Rem}[\thmcount]{Remark}
\newtheorem{Defn}[\thmcount]{Definition}
\newtheorem{Ex}[\thmcount]{Example}
\newtheorem{Asu}[\thmcount]{Assumption}
\newtheorem{Sol}[\thmcount]{Solution}
\newtheorem*{Thmx}{Theorem}
\newtheorem*{Corx}{Corollary}
\newtheorem*{Lemx}{Lemma}
\newtheorem*{Propx}{Proposition}
\newtheorem*{Remx}{Remark}
\newtheorem*{Defnx}{Definition}
\newtheorem*{Exx}{Example}
\newtheorem*{Asux}{Assumption}
\newtheorem*{Solx}{Solution}

\newtheorem{theorem}{Theorem}
\makeatletter
\newcommand{\settheoremtag}[1]{
  \let\oldthetheorem\thetheorem
  \renewcommand{\thetheorem}{#1}
  \g@addto@macro\endtheorem{
    \addtocounter{theorem}{-1}
    \global\let\thetheorem\oldthetheorem}
  }
\makeatother

\newcommand \eq[1]{\begin{equation} #1 \end{equation}}
\newcommand \eqx[1]{\begin{equation*}  #1 \end{equation*}}
\newcommand \al[1]{\begin{align} #1 \end{align}}
\newcommand \alx[1]{\begin{align*}  #1 \end{align*}}
\renewcommand \sp[1]{\begin{equation} \begin{split} #1 \end{split} \end{equation}}
\newcommand \spx[1]{\begin{equation*} \begin{split} #1 \end{split} \end{equation*}}
\newcommand \en[1]{\begin{enumerate}  #1 \end{enumerate}}

\newcommand{\thm}[2]{\begin{Thm} \label{#1} #2 \end{Thm}}

\newcommand{\lem}[2]{\begin{Lem} \label{#1} #2 \end{Lem}}

\newcommand{\prop}[2]{\begin{Prop} \label{#1} #2 \end{Prop}}

\newcommand{\defn}[2]{\begin{Defn} \label{#1} #2 \end{Defn}}

\newcommand{\pr}[1]{\begin{proof} #1 \end{proof}}

\newcounter{comcount}

\renewcommand{\hline}{\vbox{\hrule width\textwidth height 1pt}\smallskip}

       \newcommand{\be}{\beta}         \newcommand{\e}{\varepsilon}
         \newcommand{\Om}{\Omega}        \newcommand{\de}{\delta}
        \newcommand{\la}{\lambda}

\newcommand{\CC}{\mathbb{C}}

\newcommand{\NN}{\mathbb{N}}

  \newcommand{\qq}{\mathcal{Q}}
\newcommand{\RR}{\mathbb{R}}  
  
 \newcommand{\Tt}{\mathbf{T}}

\newcommand{\ZZ}{\mathbb{Z}}  

\newcommand{\supp}{\mathrm{supp}}
\newcommand{\8}{\infty}

\newcommand{\Rd}{{\RR^d}}

\renewcommand{\rm}[1]{\mathrm{#1}}
\newcommand{\wt}[1]{\widetilde{#1}}

\newcommand{\abs}[1]{\left| #1 \right|}
\newcommand{\set}[1]{\left\{ #1 \right\}}
\newcommand{\norm}[1]{\left\| #1 \right\|}

\newcommand{\eee}[1]{\left( #1 \right)}

\newcommand{\tti}[1]{T_t^{[#1]}}
\newcommand{\tsi}[1]{T_s^{[#1]}}
\newcommand{\s}[2]{%
    \IfEqCase{#2}{%
        {1}{#1^*}%
        {2}{#1^{**}}%
        {3}{#1^{***}}
    }
}
\newcommand{\dxj}{\partial_{x_j}}

\newcommand{\RL}[2]{R^{#1}_{#2}}
\newcommand{\RLQ}{\RL{j}{d_Q, loc}}
\newcommand{\RLQQ}{\RL{j}{d_Q, glob}}

\newcommand{\RDelta}[2]{\widetilde{R}^{#1}_{#2}}

\newcommand{\bx}{\mathbbm{x}}
\newcommand{\by}{\mathbbm{y}}
\renewcommand{\Tt}{T_t}
\newcommand{\Ht}{H_t}
\newcommand{\dtt}{\frac{dt}{\sqrt{t}}}
\newcommand{\dttt}{\frac{dt}{t}}
\newcommand{\dtttt}{\frac{dt}{t^2}}

\title{Riesz transform characterizations\\ for multidimensional Hardy spaces}

\author[ Edyta Kania-Strojec]{ Edyta Kania-Strojec }
 \address{
 Edyta Kania-Strojec\newline
 \indent Instytut Matematyczny, Uniwersytet Wroc\l awski \newline
 \indent pl. Grunwaldzki 2/4, 50-384 Wroc\l aw, Poland }
 \email{edyta.kania@uwr.edu.pl }

\author[ Marcin Preisner ]{ Marcin Preisner }
 \address{
 Marcin Preisner \newline
 \indent Instytut Matematyczny, Uniwersytet Wroc\l awski \newline
 \indent pl. Grunwaldzki 2/4, 50-384 Wroc\l aw, Poland }
 \email{marcin.preisner@uwr.edu.pl }

\subjclass[2010]{}
\thanks{ The authors are supported by the grant No. 2017/25/B/ST1/00599 from National Science Centre (Narodowe Centrum Nauki), Poland. }

\newcommand{\Har}{H^1_L(X)}
\newcommand{\sumQ}{\sum_{Q\in\qq}}
\newcommand{\sleq}{\lesssim}
\newcommand{\sgeq}{\gtrsim}

\begin{document}
\maketitle

\begin{abstract}
We study Hardy space $H^1_L(X)$ related to a self-adjoint operator $L$ defined on an Euclidean subspace $X$ of $\Rd$. We continue study from \cite{KPP}, where, under certain assumptions on the heat semigroup $\exp(-tL)$, the atomic characterization of local type for $H^1_L(X)$ was proved.

In this paper we provide additional assumptions that lead to another characterization of $H^1_L(X)$ by the Riesz transforms related to $L$. 
As an application, we prove the Riesz transform characterization of $H^1_L(X)$ for multidimensional Bessel and Laguerre operators, and the Dirichlet Laplacian on $\RR^d_+$.
\end{abstract}

\newcommand{\HarRd}{H^1(\RR^d)}

\section{Introduction and statement of results}

\subsection{Introduction.}
Let $H_t = \exp(t\Delta)$ be the heat semigroup on $\Rd$, i.e. $H_t f(x) = \int_{\Rd} H_t(x-y) f(y)  \, dy$ and
\eq{\label{heat_kernel}
H_t(x-y) = (4\pi t)^{-d/2} \exp\eee{-\frac{|x-y|^2}{4t}}, \quad x,y\in \Rd, \ t>0.
}
The classical Hardy space $\HarRd$ can be defined by the maximal operator related to the operators $H_t$ and plays an important role in harmonic analysis. We say that a function $f\in L^1(\Rd)$ is in $ \HarRd$ if and only if
\eqx{\label{eq_class}
\norm{f}_{\HarRd} := \norm{ \sup_{t>0} \abs{H_t f(\cdot)}}_{L^1(\RR^d)}<\8.
}
There are many equivalent definitions of $\HarRd$ related to various objects in harmonic analysis. The interested reader is referred to \cite{Stein} and references therein. Let us recall that the Riesz transforms $\widetilde{R}_j = \dxj  (-\Delta)^{-1/2}$, $j=1,...,d$, are given by
\spx{
\widetilde{R}_jf(x) & = C_d \lim_{\e\to 0} \int_{|x-y|>\e} \frac{x_j-y_j}{|x-y|^{d+1}} f(y) \, dy,
}
where $x= (x_1,...,x_d)\in \Rd$. One of the classical results states that one can give equivalent definition of $\HarRd$ in terms of the Riesz transforms, c.f. \cite{Fefferman_Stein}. More precisely a function $f$ belongs to $H^1(\Rd)$ if and only if all the functions: $f,\widetilde{R}_1f,...,\widetilde{R}_d f$ belong to $L^1(\Rd)$ and
\eq{\label{eq_class_Riesz}
\norm{f}_{\HarRd}\simeq \norm{f}_{L^1(\Rd)} + \sum_{j=1}^d \norm{\wt{R}_jf}_{L^1(\Rd)}.
}

On the other hand, a function $f$ in $\HarRd$ can be decomposed as an infinite linear combination of simple functions called atoms, see \cite{Coifman_Studia} and \cite{Latter_Studia}. More precisely, for a function $f\in \HarRd$ we can write
\eq{\label{class_decomp}
f = \sum_{k=1}^\8 \la_k a_k,
}
where $\sum_k |\la_k| <\8$ and $a_k$ are {\it atoms}, i.e. there exist balls $B_k$ in $\Rd$ such that:
\eq{\label{class_atoms}
\supp \, a_k \subseteq B_k, \qquad \norm{a_k}_\8 \leq |B_k|^{-1}, \qquad \int_{B_k} a_k (x) dx = 0.
}
Here $|B_k|$ is the Lebesgue measure of the ball $B_k$. For more properties of $\HarRd$ we refer the reader to \cite{Stein} and references therein.

One can consider $\HarRd$ as related to the classical Laplacian $\Delta$ on $\Rd$, since many possible definitions of $\HarRd$ are given in terms of $\Delta$. Since the 60's many researchers considered the Hardy spaces $H^1_L(X)$ related to various self-adjoint operators $L$ on some metric-measure spaces $X$, see e.g. \cite{ Uchiyama, Latter_Studia, Song_Yan_2016, CoifmanWeiss_BullAMS, Fefferman_Stein, Bernicot, BDT_d'Analyse, DP_Annali, DZ_Potential_2014, Dziubanski_Constructive, DH_Dunkl,Hofmann_Memoirs , Goldberg_Duke}. A natural question in this theory is the following: can we have decompositions of the type \eqref{class_decomp} for $f\in H^1_L(X)$? Also,  whether the equivalence similar to \eqref{eq_class_Riesz} holds or not? It appears that now we have many general results concerning atomic decompositions for $H^1_L(X)$, see e.g. \cite{Hofmann_Memoirs, DP_Annali, PSY, Uchiyama}. However, the characterization of $H^1_L(X)$ in terms of the Riesz transforms is not known in such generality.

In the present paper we shall continue study in the context introduced in \cite{KPP}. Recall, that in \cite{KPP} we consider a space $X \subseteq \Rd$ and a nonnegative self-adjoint operator $L$ on  $L^2(X)$. The semigroup $\exp(-tL)$ satisfy the upper Gaussian estimates and, roughly speaking, the kernel $T_t(x,y)$ of $\exp(-tL)$ is similar to $H_t(x-y)$ for local times and $T_t(x,y)$ decays faster for global times, where the scale of time is adjusted to some covering $\qq = \{Q_j\}_{j\in \NN}$ of $X$. For a precise statement of these assumptions see \cite{KPP} or Section \ref{sec_assumptions} below. The main issue considered in \cite{KPP} was the characterization of $H_L^1(X)$ in terms of the atomic decompositions. It was proved there  that in this context one has atoms for $H^1_L(X)$ that are either classical atoms (as in \eqref{class_atoms}) or atoms of the form $a(x) = |Q|^{-1}\mathbbm{1}_Q(x)$, $Q\in \qq$. The latter atoms are called ''local atoms'', c.f.\ \cite{Goldberg_Duke}.

Our goal here is to  characterise $\Har$ by the Riesz transforms $D_j L^{-1/2}$, $j=1,...,d$, where $D_j = \partial_{x_j} +V_j$ is a derivative adapted to $L$. To this end we add additional assumptions for the kernels: $\partial_{x_j}T_t(x,y)$, $V_j(x) T_t(x,y)$. Using this we show a result similar to \eqref{eq_class_Riesz}, i.e.\ the Hardy space $H^1_L(X)$ is characterized by appropriate Riesz transforms. For other results concerning this question, see e.g. \cite{BDG_Tohoku, BDT_d'Analyse, Preisner_JAT, Fefferman_Stein, MR3639532, Dziubanski_Bessel_Dunkl,Hejna_MathNachr,  DP_Potential, MR2452309}.

Our main motivation here is to give an uniform approach that will work in different contexts  to study operators such as: multidimensional Bessel and Laguerre operators, and the Dirichlet Laplacian on $\RR^d_+$. In the last and most technical section we verify that our assumptions are indeed satisfied for these examples.

\subsection{Assumptions}\label{sec_assumptions}
In this section we state assumptions that will be used throughout the paper. Let $X\subseteq \Rd$ be a space that is a product of: finite intervals, half-lines, or lines equipped with the Lebesgue measure, i.e.\ $X= (a_1,b_1)\times...\times(a_d,b_d)$, where $a_j \in [-\8, \8)$ and $b_j \in (-\8,\8]$. We shall study a non-negative self-adjoint operator $L$ that is densely defined on $L^2(X)$. The semigroup generated by $-L$ will be denoted by $T_t = \exp(-tL)$ and we further assume that there exists an integral kernel $T_t(x,y)$, such that for $f \in L^p(X)$, $1\leq p \leq \8$, we have
$$T_t f(x) = \int_X T_t(x,y) f(y)\, dy, \qquad \text{a.e. }x\in X.$$

The Hardy space $\Har$ related to $L$ is defined in terms of the maximal operator related to $T_t$, namely
\eqx{ \label{max-hardy}
	\Har = \set{ f\in L^1(X) \ : \ \norm{f}_{\Har} := \norm{\sup_{t>0} \abs{ \Tt f}}_{L^1(X)}  < \8}.
	}

In this paper we shall study the spaces $H^1_L(X)$  that will be related to some coverings $\qq = \set{Q_k \ : \ k\in \NN}$ of $X$, where $Q_k$ are  cuboids. We assume that $\qq$ is {\it an admissible covering} in the sense of Definition \ref{def_covering} below. Let $d_Q$ be the diameter of $Q$ and denote by $Q^*$ a slight enlargement of $Q$, see the comments after Definition \ref{def_covering} below. Following \cite{KPP} we assume that there exists $\gamma \in (0,1/3)$ and $C,c>0$, such that $\Tt(x,y)$ satisfies:
\begin{align}
\label{a0} \tag{$A_0$}
&0\leq \Tt(x,y) \leq C t^{-d/2} \exp\left(-\frac{|x-y|^2}{c t}\right), && x,y \in X, t>0,\\
 \label{a1} \tag{$A_1$}
&\sup_{y \in \s{Q}{2}} \int_{(\s{Q}{3})^c} \sup_{t > 0} t^\delta \Tt(x,y) \, dx \leq C d_Q^{2\delta}, && \delta \in [0,\gamma),  Q\in \qq,
\\
\label{a2} \tag{$A_2$}
&\sup_{y \in \s{Q}{2}} \int_{\s{Q}{3}} \sup_{t \leq  d_Q^2 } t^{-\delta} \abs{\Tt(x,y) - H_t(x-y)} \,dx \leq C d_Q^{-2\delta}, && \delta \in [0,\gamma), Q\in \qq.
\end{align}
In \cite{KPP} the authors studied $H^1_L(X)$ for operators satisfying \eqref{a0}--\eqref{a2}. It was proved that $H^1_L(X)$ can be characterized by atomic decompositions with local atoms of the form $|Q|^{-1} \mathbbm{1}_Q$, where $Q\in \qq$, see  \cite[Thm.\ A]{KPP} and Theorem \ref{thm_decomp} below.

In the present paper we shall study the Riesz transform characterization of $\Har$, when $L$ satisfies \eqref{a0}--\eqref{a2} and the following assumptions that are inspired by certain known examples like: Bessel, Laguerre, or Schr\"odinger operators. On $L^2(X)$ consider the operators $R_j$ formally given by:
 $$R_j = (\partial_{x_j}+V_j) L^{-1/2}, \qquad j=1,...,d,$$
 where $\partial_{x_j}$ is the standard derivative and $V_j$ is a function that depends only on $x_j$. Suppose that $T_t(x,y)$ satisfy:
 \al{\label{a3} \tag{$A_3$} 
&\sup_{y\in Q^{**}} \int_{(Q^{***})^c} \int_0^{d_Q^2} \abs{\partial_{x_j} T_t(x,y)} \, \frac{dt}{\sqrt{t}} \, dx  \leq C, && Q\in\qq, j=1,...,d,
\\
\label{a4} \tag{$A_4$}
&\sup_{y\in Q^{**}} \int_X \int_{d_Q^2}^\8 \abs{\dxj  T_t(x,y)} \, \frac{dt}{\sqrt{t}} \, dx \leq C, && Q\in \qq,  j=1,...,d,
\\
\label{a5} \tag{$A_5$}
	&\sup_{y\in Q^{**}} \int_{Q^{***}} \int_0^{d_Q^2}  \abs{\dxj \eee{T_t(x,y) -  H_t(x-y)} } \, \frac{dt}{\sqrt{t}} \, dx\leq C, && Q\in\qq,  j=1,..., d,
	\\
	\label{a6} \tag{$A_6$} &\sup_{y\in X} \int_X \int_0^\8  \abs{V_j(x)} T_t(x,y) \, \dtt \, dx \leq C, && j=1,...,d.
	}
For  $j=1,...,d$ define the kernels
 \newcommand{\inte}{\int_\e^\e}
\eq{\label{eq_R_ker}
R_j(x,y) := \pi^{-1/2} \int_0^\8  \eee{\dxj + V_j(x_j)}  \Tt(x,y) \dtt.
}
Notice that our assumptions guarantee that the integral above exists for a.e. $(x,y)$. The operators $R_j$ are defined as follows:
\eqx{\label{stowarzyszenie}
R_j f(x) = \lim_{\e\to 0} \int_{|x-y|>\e} R_j(x,y) f(y) \, dy, \qquad x\in X.
}
We always assume that $R_j$ are bounded on $L^2(X)$.


\subsection{Results} \label{ssec_1.3}
Our first main result is the following theorem, that describes the Hardy space $H^1_L(X)$ in terms of the Riesz transforms.
\settheoremtag{A}
\begin{theorem} \label{thma} Assume that there is an operator $L$ and an admissible covering $\qq$ as in Section \ref{sec_assumptions}. In particular, we assume that \eqref{a0}--\eqref{a6} are satisfied. Then $f\in \Har$ if and only if $f,R_1f,...,R_d f\in L^1(X)$. Moreover, there exists a constant $C>0$ such that
\eqx{
C^{-1} \norm{f}_{\Har} \leq \norm{f}_{L^1(X)} + \sum_{j=1}^d \norm{R_j f}_{L^1(X)} \leq C  \norm{f}_{\Har}.
}
\end{theorem}

The proof of Theorem \ref{thma} is given in Section \ref{proof-D} below and it is based on known techniques. The main idea is to compare (locally) $R_j$ with the classical Riesz transforms $\wt{R}_j =\dxj  (-\Delta)^{1/2}$ and use additional decay as $t\to \8$.  

One of our main motivations is to study product cases. Assume that for $i=1,...,N$ we have operators $L_i$ satisfying the assumptions as in Section \ref{sec_assumptions}. In particular, $L_i$ is associated with the semigroup $\tti{i}$ that has a kernel $\tti{i}(x_i,y_i), \ x_i,y_i \in X_i$. Then we can define
\eq{\label{prod_X}
X = \prod_{i=1}^N X_i \subseteq \prod_{i=1}^N \RR^{d_i} = \Rd
}
and
\eq{\label{prod_L}
L  = L_1 + ... + L_N,
}
such that  each $L_i$ acts only on the variable $x_i \in X_i$. For more precise description see Section \ref{sec_products} below. The following theorem gives the Riesz transform characterization for $H^1_{L}(X)$ in the product case. 

\settheoremtag{B}
\begin{theorem}\label{thmb}
Let $X$ and  $L$ be as in \eqref{prod_X}--\eqref{prod_L} and assume that for each $i=1,...,N$ the semigroup kernel $\tti{i}(x_i, y_i)$ together with an admissible covering $\qq_i$ of $X_i$ satisfy the conditions  \eqref{a0}--\eqref{a6}. Then $f\in H^1_{L}(X)$ if and only if $f, R_1 f, ..., R_d f \in L^1(\Rd)$. Moreover,
$$C^{-1} \norm{f}_{H^1_{L}(X)} \leq \norm{f}_{L^1(\Rd)} + \sum_{j=1}^d \norm{R_j f}_{L^1(\Rd)} \leq C \norm{f}_{H^1_{L}(X)}. $$
\end{theorem}

The proof of Theorem \ref{thmb} is given in Section \ref{proof-E} below. We shall use \cite[Thm.~B]{KPP}, where we proved that assuming \eqref{a0}--\eqref{a2} for $\tti{i}(x_i,y_i)$ and $\qq_i$ we can define an admissible covering $\qq_1 \boxtimes ... \boxtimes Q_N$ that describes $H^1_L(X)$ for $L=L_1+...+L_N$, see  \cite[Def.~1.5]{KPP}.

As an example of applications of Theorem \ref{thmb} we study certain multidimensional Bessel and Laguerre operators. Thanks to Theorem \ref{thmb} it is enough to verify \eqref{a0}--\eqref{a6} only in the one-dimensional case. Then, the Riesz transform characterization for $H^1_L(X)$ for the multidimensional case (when $L$ is the sum of Bessel or Laguerre operators) follows from Theorem \ref{thmb}. Furthermore, a similar argument (see Section \ref{sec3.1}) allows us to study the Dirichlet Laplacian on the half-space $\RR^d_+$. Below we briefly recall the operators that we work with and state the results.

\newcommand{\lbl}{L_B^{[\be]}}
\newcommand{\lla}{L_L^{[\be]}}
{\bf Bessel operator.} Let $X=(0,\8)^d$. For $\be=(\beta_1,...,\be_d)$ assume $\be_i>0$, $i=1,...,d$, and consider the multidimensional Bessel operator
\eq{\label{L_Bess}
\lbl=- \sum_{i=1}^d \eee{\frac{d^2}{dx_i} - \frac{\be_i^2-\be_i}{x_i^2}}, \qquad x_1,...,x_d>0.
}
More precisely, by $\lbl$ we shall denote a proper self-adjoint operator defined on $L^2(X)$, see e.g. \cite{MR2284989}. Harmonic analysis related to $\lbl$ was studied in e.g. \cite{BDT_d'Analyse, Betancor_Buraczewski, Betancor_Castro_Nowak, MR2600427, MR1808644, MR2284989}. In \cite{BDT_d'Analyse} the authors describe the Hardy space related to $\lbl$ for $d=1$ in terms of either atomic decompositions or Riesz transforms 
$$R_j = \eee{\partial_{x_j}-\frac{\beta_j}{x_j}} \eee{L_B^{[\be]}}^{-1/2}, \qquad j=1,...,d.$$
Denote
\eq{\label{Bess_cov}
\qq_B = \set{[2^n,2^{n+1}] \, : n\in \ZZ}.
}
Then $\qq_B$ is an admissible covering for $(0,\8)$ and for $d>1$ we have the admissible coverings $\qq_B \boxtimes ... \boxtimes \qq_B$ defined in \cite[Def.\ 1.5]{KPP}. The following theorem follows directly from \cite[Prop.\ 4.3]{KPP}, Theorem \ref{thmb} and Proposition \ref{riesz_bessel} below.

\settheoremtag{C}
\begin{theorem}\label{coro_bess-riesz}
Let $d\geq 1$, $\be_1,...,\be_d>0$ and $\lbl$ be the multidimensional Bessel operator, see \eqref{L_Bess}. Then, $f\in H^1_{\lbl}\eee{(0,\8)^d}$ if and only if $f,R_1 f, ..., R_d f \in L^1\eee{(0,\8)^d}$. Moreover, the associated norms are comparable, i.e.
$$\norm{f}_{H^1_{\lbl}\eee{(0,\8)^d}} \simeq \norm{f}_{L^1((0,\8)^d)} + \sum_{j=1}^d \norm{R_jf}_{L^1((0,\8)^d)}.$$
\end{theorem}

{\bf Laguerre operator.}
Let $\be=(\be_1,...,\be_d)$, where  $\be_i>0$, $i=1,...,d$, and denote the multidimensional Laguerre operator
\eq{\label{L_Lagg}
\lla=- \sum_{i=1}^d \eee{\frac{d^2}{dx_i^2}-x_i^2-\frac{\be_i^2-\be_i}{x_i^2}}, \qquad x_1,...,x_d>0.
}
Set $X=(0,\8)^d$. By $\lla$ we shall denote a known self-adjoint operator on $L^2(X)$, see  e.g. \cite{Nowak_Stempak_Advances} . In \cite{BDG_Tohoku, MR2786703,  NowakStempak_hermite, MR2381163, Nowak_Stempak_Advances} we find some studies on harmonic analysis related to $\lla$. In particular the authors of \cite{BDG_Tohoku} proves the atomic decomposition theorem for the Hardy space related to $\lla$ in the one-dimensional case. For $d=1$ we have the following admissible covering of $(0,\8)$,
\sp{\label{Lag_cov}
\qq_{L}=&\set{[2^n+(k-1)2^{-n},2^n+k2^{-n}]\colon k=1,\ldots,2^{2n}; n\in\NN }\\
&\cup\set{[2^{-n},2^{-n+1}]\colon n\in \NN_+}.}
and, using this covering, we produce $\qq_L \boxtimes ... \boxtimes \qq_L$ for $d>1$, see \cite[Def.\ 1.5]{KPP}. Combining \cite[Prop.\ 4.5]{KPP}, Prop.\ \ref{a_lag_riesz} below, and Theorem \ref{thmb} we arrive at the following characterization of $H^1_{L_L^{[\be]}}((0,\8)^d)$ in terms of the Riesz transforms
$$R_j = \eee{\partial_{x_j}+x_j-\frac{\be_j}{x_j}} \eee{L_L^{[\be]}}^{-1/2}.$$
\settheoremtag{D}
\begin{theorem}\label{coro_lag}
Let $d\geq 1$, $\be_1,...,\be_d>0$ and $\lla$ be the multidimensional Laguerre operator, c.f.\ \eqref{L_Lagg}. Then, $f\in H^1_{\lla}\eee{(0,\8)^d}$ if and only if $f,R_1 f, ..., R_d f \in L^1\eee{(0,\8)^d}$. Moreover, the associated norms are comparable, i.e.
$$\norm{f}_{H^1_{\lla}\eee{(0,\8)^d}} \simeq \norm{f}_{L^1((0,\8)^d)} + \sum_{j=1}^d \norm{R_jf}_{L^1((0,\8)^d)}.$$
\end{theorem}

{\bf Dirichlet Laplacian on $\RR^d_+$.}
As a third example let us consider the  Dirichlet Laplacian on the half-space $\RR^d_+$. To be more precise for $d\geq 1$ we consider
$$X=R^d_+ =\set{x=(x_1,...,x_d)\in \Rd \ : x_d>0}$$ and the Laplacian $L_D = -\Delta$ with the Dirichlet boundary condition at $x_d=0$. The semigroup generated by $-L_D$ is given by $T_t f(x) = \int_X T_t(x,y) f(y) \, dy$,
$$T_t(x,y) = H_t(x-y) - H_t(x-\wt{y}), \quad x,y\in X$$
where $\wt{y} = (y_1,...,y_{d-1}, -y_d)$ and $H_t$ is as in \eqref{heat_kernel}. In this case the appropriate covering $\qq_{D}$ is a covering that consists of cubes such that a cube $Q\in \qq_{D}$ has the diameter comparable to the distance from the boundary $x_d=0$ (in the case $d=1$ one can just take dyadic covering of $(0,\8)$). The Hardy space for this operator was studied in e.g. \cite{Auscher_Russ, Chang_Krantz_Stein, PSY}. It is not hard to check that  \eqref{a0}--\eqref{a5} are satisfied (see Section \ref{sec:ex-dir-lap}) and we obtain the following characterization of $H^1_{L_D}(\RR^d_+)$ by means of the Riesz transforms $R_j = \partial_{x_j} L_D^{-1/2}$, $j=1,...,d$.
\settheoremtag{E}
\begin{theorem}\label{coro_dir}
Let $d\geq 1$ and $L_D$ be the Laplacian on $\RR^d_+$ with the Dirichlet boundary condition at $x_d=0$. Then, $f\in H^1_{L_D}\eee{\RR^d_+}$ if and only if $f,R_1 f, ..., R_d f \in L^1\eee{\RR^d_+}$. Moreover, the associated norms are comparable, i.e.
$$\norm{f}_{H^1_{L_D}(\RR^d_+)} \simeq \norm{f}_{L^1(\RR^d_+)} + \sum_{j=1}^d \norm{R_jf}_{L^1(\RR^d_+)}.$$
\end{theorem}

{\bf Organization of the paper.} In Section \ref{sec2} we recall some known facts and prove preliminary estimates. Theorems \ref{thma} and \ref{thmb} are proved in Section \ref{sec3}. In Section \ref{sec3.1} we state and prove a modification of Theorem \ref{thmb} that will be needed in the proof of Theorem \ref{coro_dir}. Propositions \ref{riesz_bessel} and \ref{a_lag_riesz}, that are crucial for Theorems \ref{coro_bess-riesz} and \ref{coro_lag}, are stated and proved in Section~\ref{sec4}. At the end of Section \ref{sec4} we briefly discuss the Dirichlet Laplacian and prove Theorem \ref{coro_dir}. We shall use a standard convention that   $C$ and $c$ at each occurrence denote
some positive constants independent of relevant quantities (depending on the context). We will write $A \sleq B$ for $A \leq CB$ and $A \simeq B$ for $A\sleq B \sleq A$.

\section{Preliminaries}\label{sec2}

\subsection{Admissible coverings}

Let $X \subseteq \Rd$ be as in Section\ \ref{sec_assumptions}. For $z = (z_1, ..., z_d)\in X$ and $r_1,...,r_d>0$ we denote the closed cuboid
$$Q(z,r_1,...,r_d) = \set{x\in X \ : \ |x_i-z_i|\leq r_i \text{ for } i=1,...,d},$$
and the cube
$Q(z,r) = Q(z,r,...,r).$ The following definition  will be used throughout the paper, c.f. \cite[Def.\ 1.2]{KPP}.
\defn{def_covering}{
Let $\qq$ be a  set of cuboids in $X\subseteq \Rd$. We call $\qq$ {\it an admissible covering} if:
\en{
\item $X = \bigcup_{Q\in \qq} Q$,
\item if $Q_1, Q_2 \in \qq$ and $Q_1\neq Q_2$, then $|Q_1\cap Q_2| = 0$,
\item if $Q = Q(z,r_1,...,r_d)\in \qq$, then $r_i\simeq r_j$ for $i,j \in \set{1,...,d}$,
\item if $Q_1, Q_2 \in \qq$ and $Q_1 \cap Q_2 \neq \emptyset$, then $ d_{Q_1} \simeq d_{Q_2}$,
\item if $Q\in \qq$, then $\mathrm{dist_{\Rd}}(Q, \Rd \setminus X) \sgeq d_Q$.
}
} 
Having an admissible covering $\qq$ and $Q=(z,r_1,...,r_d) \in \qq$, we define 
$$Q^*:= Q(z, \kappa r_1,...,\kappa r_d),$$
where $\kappa>1$ is chosen so that for $Q_1,Q_2 \in \qq$,
\eqx{\label{neighbours}
Q_1^{***}\cap Q_2^{***} \neq \emptyset \quad \iff \quad  Q_1 \cap Q_2 \neq \emptyset}
and 
\eqx{\label{two_enlargements}
\mathrm{dist}_{\Rd}(Q^{***},\Rd\setminus X)>0.
}
The family $\set{Q^{***}}_{Q\in\qq}$ is a finite covering of $X$, namely
\eq{\label{finite_covering}
\sum_{Q\in \qq} \mathbbm{1}_{Q^{***}}(x) \leq C, \qquad x\in X. 
}
Let us notice that we have a flexibility in choosing the enlargements $Q^*$, $Q^{**}$, $Q^{***}$ etc. In particular the notation in \cite{KPP} is slightly different. Recall that having admissible coverings $\qq_i$ of $X_i$, $i=1,...,N$, we can produce a natural admissible covering $\qq_1\boxtimes ... \boxtimes \qq_N$ of $X$ as in \eqref{prod_X}, see \cite[Def. 1.5]{KPP}.

\subsection{Products}\label{sec_products}

In this subsection $i$ will be always an index from $\set{1,...,N}$. Let $X_i \subseteq \RR^{d_i}$ and $L_i$ are as in Section \ref{sec_assumptions} on $L^2(\RR^{d_i})$. Set $d=d_1+...+d_N$ and let $X$ be as \eqref{prod_X}. Now, we shall explain the precise meaning of \eqref{prod_L}. 
Slightly abusing the notation we keep the symbol $L_i$ for the operator
\eqx{  \underbrace{I\otimes ... \otimes I }_{i-1 \text{ times}}\otimes L_i \otimes  \underbrace{I \otimes ... \otimes I}_{N-i \text{ times}} }
 on $L^2(X)$, where $I$ denotes the identity operator on the corresponding subspace, and we define
\eqx{\label{sum_of_operators} Lf(x) = L_1f(x) + \ldots + L_Nf(x), \qquad x=(x_1,\ldots,x_N)\in X.}
Since the operators $L_i$ are self-adjoint, the operator $L$ is well defined and essentially self-adjoint, see e.g. \cite[Thm.\ 7.23]{MR2953553}.

Recall that the semigroups $\tti{i} = \exp\eee{-tL_i}$ on $X_i$ have the kernels $\tti{i}(x_i,y_i)$, $x_i,y_i \in X_i, t>0$, so that the semigroup $T_t = \exp(-tL)$ is related to the kernel
$$T_t(x,y) = T_t^{[1]}(x_1,y_1) \cdot ... \cdot  T_t^{[N]}(x_N,y_N).$$

\subsection{Local atomic Hardy spaces}\label{sec_atomic}

For an admissible covering $\qq$ of $X\subseteq \Rd$ (see Definition \ref{def_covering}) we shall define {\it  the local atomic Hardy space} $H^1_{at}(\qq)$ related to  $\qq$ as follows.
\defn{def_Qatoms}
{
A function $a \, : \, X \to \CC$ is called a $\qq-atom$ if either:
\begin{itemize}
\item[(i)]  there is  $Q\in\qq$ and a cube $K\subset Q^{**}$, such that:  $$\supp \, a \subseteq K, \ \ \norm{a}_\8 \leq |K|^{-1},  \ \ \int a(x)\, dx = 0;$$
\end{itemize}
or
\begin{itemize}
\item[(ii)] there exists $Q\in\qq$ such that $$a(x) = |Q|^{-1}\mathbbm{1}_Q(x).$$
\end{itemize}
}

Then, the atomic space $H^1_{at}(\qq)$, is defined in a standard way. Namely, we say that a function $f$ is in $H^1_{at}(\qq)$ if $f = \sum_k \la_k a_k$ with $\qq$-atoms $a_k$ and $\sum_k |\la_k| <\8$. Moreover, the norm of $H^1_{at}(\qq)$ is given by
$$\norm{f}_{H^1_{at}(\qq)} = \inf \sum_k \abs{\la_k},$$
where the infimum is taken over all possible representations of $f = \sum_k \la_k a_k$ as above. A standard argument shows that $H^1_{at}(\qq)$ is a Banach subspace of $L^1(X)$.

Here we state the atomic decomposition result that follows from \cite[Thm.~A]{KPP}. This will be needed later on in the proof of Theorem \ref{thma}.
\thm{thm_decomp}{
Assume that for $L, T_t$, and an admissible covering $\qq$ the assumptions \eqref{a0}--\eqref{a2} are satisfied. Then $H^1_L(X) = H^1_{at}(\qq)$ and the corresponding norms are equivalent.
}

\newcommand{\Rtau}{\wt{R}^j_{\tau, loc}}
\newcommand{\Rtauu}{\wt{R}^j_{\tau, glob}}
\subsection{Classical local Hardy spaces.} \label{ssec24}

In this section we recall briefly some theory related to the classical local Hardy spaces on $\Rd$, c.f. \cite{Goldberg_Duke, Stein}. In particular, we shall present the relation between classical local Hardy spaces and local Riesz transforms in Proposition \ref{prop_local_riesz}.

Recall that the kernel of the Riesz transform $\wt{R}_j = \dxj (-\Delta)^{-1/2}$ can be given by $\wt{R}_j(x,y) = \pi^{-1/2} \int_0^\8 \dxj H_t(x-y)\, \dtt$ and for $\tau>0$ denote
$$\Rtau(x,y) = \pi^{-1/2} \int_0^{\tau^2} \dxj  \Ht(x-y) \frac{dt}{\sqrt{t}}, \qquad \Rtauu(x,y) = \pi^{-1/2} \int_{\tau^2}^\8 \dxj  \Ht(x-y) \frac{dt}{\sqrt{t}} .$$
It is well known that these kernels are related (in the principal value sense) with the operators $\Rtau$ and $\Rtauu$ that are well-defined and  bounded on $L^2(\Rd)$ (uniformly in $\tau>0$). In what follows we shall need the following version of the characterization of the local Hardy spaces.


\prop{prop_local_riesz}{
There exist $C_1,C_2>0$ that does not depend on $\tau>0$ such that:
 \en{
 \item If $a$ is either a classical atom or local atom of the form $a=|Q|^{-1}\mathbbm{1}_{Q}$, where $Q=Q(z,r_1,...,r_d)$, $r_1\simeq ...\simeq r_d \simeq \tau$, we have
 	$$\norm{a}_{L^1(\Rd)} + \sum_{j=1}^d \norm{\Rtau a}_{L^1(\Rd)} \leq C_1.$$
 	
 \item Assume that $\supp f \subseteq Q^*$, where $Q=Q(z,r_1,...,r_d)$, $r_1\simeq ...\simeq r_d \simeq \tau$, and
 	$$M:=\norm{f}_{L^1(Q^*)} + \sum_{j=1}^d \norm{\Rtau f}_{L^1(Q^{**})}<\8.$$
Then there exist sequences $\{\la_k\}_k$ and $\{a_k\}_k$, such that: $f = \sum_k \la_k a_k$, $\sum_k \abs{\la_k} \leq C_2 M$, and $a_k$ are either the classical atoms supported in a cube $K \subseteq Q^{**}$ or $a_k = |Q|^{-1}\mathbbm{1}_{Q}$.
 }
 }

\pr{[Sketch of the proof] 
This fact is well known and has quite standard proof. For the convenience of the reader we provide a sketch of the proof. Notice that
$$\Rtau (x,y) = c_d \frac{x_j-y_j}{|x-y|^{d+1}} \psi\eee{\frac{|x-y|}{\tau}},
$$
where $\psi$ is smooth on $[0,\8)$, $\psi(0) = c_d'$ and $\psi(s) \simeq e^{-s^2}$ as $s \to \8$.

Part {\bf 1.} follows by a standard Calder\'on-Zygmund argument. The main idea is to use the $L^2$-estimate on $Q(x_0, 2\tau)$ and the estimate $\Rtau (x,y) \leq \tau|x-y|^{-d-1}$ for $y\in Q(x_0,\tau)$ and $x\not\in Q(x_0,2\tau)$. 

In order to prove {\bf 2.} define $\la_0 = \int f$ and let
$$g(x) = f(x) - \la_0 |Q|^{-1} \mathbbm{1}_Q(x).$$
Then $a_0(x) = |Q|^{-1} \mathbbm{1}_Q(x)$ is one of our atoms, $|\la_0| \leq M$, $\supp \, g \subseteq Q^*$ and $\int g = 0$. By standard computations one may check that
$$\norm{g}_{L^1(\Rd)} + \sum_{j=1}^d \norm{\wt{R}_j g}_{L^1(\Rd)} \sleq M.$$
Using the classical characterization of $H^1(\Rd)$ by means of the Riesz transforms, see \eqref{eq_class_Riesz}, we obtain
$$g(x) = \sum_{k=1}^\8 \la_k a_k(x),$$
where $a_k$ are classical atoms on $\Rd$ and 
$$\sum_{k=1}^\8 |\la_k| \sleq M.$$
Then
$$f(x) = \sum_{k=0}^\8 \la_k a_k(x), \qquad \sum_{k=0}^\8 |\la_k| \sleq M.$$
This may look that we are done, but notice that we also want to have atoms $a_k$ supported in $Q^{**}$ (not anywhere in $\Rd$). This can be done by a standard procedure, for details see e.g. \cite[Thm.\ 2.2(b)]{Kania1}. Let us notice, that here we make use of the property {\bf 5.} of Definition~\ref{def_covering}, i.e.\ we enlarge $Q$ in $\Rd$, but we want to have atoms supported in $Q^{**}$ that is still in $X$.
}

\subsection{Partition of unity}\label{sec_partition}

In what follows we shall decompose functions using an admissible covering $\qq$ of $X\subseteq \Rd$ see Definition \ref{def_covering}. One can find functions $\psi_Q\in C^1(X)$ such that:
\eqx{\label{partition}
0 \leq \psi_Q (x) \leq \mathbbm{1}_{Q^{*}}(x), \quad \norm{\psi_Q'}_\8 \leq C d_Q^{-1}, \quad \sum_{Q\in \qq} \psi_Q(x) = \mathbbm{1}_X(x).
}
The family $\set{\psi_Q}_{Q\in\qq}$ will be called  {\it a partition of unity} related to $\qq$.

\subsection{Auxiliary estimates.}
In what follows we shall use a slight generalization of \eqref{a2}--\eqref{a5} that we state below for further references. 

\lem{constants}{Assume that $T_t$ together with admissible covering $\qq$ satisfy \eqref{a0} and \eqref{a2} -- \eqref{a5}. Let $\gamma$ be as in \eqref{a2}. Then, for $c\geq 1$ there exists $C>0$ such that
\al{\label{a2prim} \tag{$A_2'$}
&\sup_{y \in \s{Q}{2}} \int_{\s{Q}{3}} \sup_{t \leq c d_Q^2 } t^{-\delta} \abs{\Tt(x,y) - H_t(x-y)} \,dx \leq C d_Q^{-2\delta}, && \delta \in [0,\gamma), Q\in \qq.\\
\label{a3prim} \tag{$A_3'$}
&\sup_{y\in Q^{**}} \int_{(Q^{***})^c} \int_0^{cd_Q^2} \abs{\partial_{x_j} T_t(x,y)} \, \frac{dt}{\sqrt{t}} \, dx  \leq C, && Q\in\qq, j=1,...,d,
\\\label{a4prim} \tag{$A_4'$}
&\sup_{y\in Q^{**}} \int_X \int_{c^{-1}d_Q^2}^\8 \abs{\dxj  T_t(x,y)} \, \frac{dt}{\sqrt{t}} \, dx \leq C, && Q\in \qq,  j=1,...,d,
\\\label{a5prim} \tag{$A_5'$}
&\sup_{y\in Q^{**}} \int_{Q^{***}} \int_0^{cd_Q^2}  \abs{\dxj \eee{T_t(x,y) -  H_t(x-y)} } \, \frac{dt}{\sqrt{t}} \, dx\leq C, && Q\in\qq,  j=1,..., d.
}}

The proof of Lemma \ref{constants} is a simple exercise that follows easily from \eqref{a0}  and \eqref{a2}--\eqref{a5}.

\subsection{Riesz transforms.}

For $\tau >0$ and $j=1,...d$ we split the kernel \eqref{eq_R_ker} as $R_j(x,y) = \RL{j}{\tau,loc}(x,y) +  \RL{j}{\tau, glob}(x,y) +  R^j_{V}(x,y) $, where
\sp{\label{loc_glob}
\RL{j}{\tau,loc}(x,y) = \pi^{-1/2} \int_0^{\tau^2} \dxj T_t(x,y)\, \dtt,\qquad x,y\in X, \\
\RL{j}{\tau, glob}(x,y) = \pi^{-1/2} \int_{\tau^2}^\8 \dxj T_t(x,y)\, \dtt,\qquad x,y\in X,\\
R^j_{V}(x,y) = \pi^{-1/2} \int_0^\8 V_j(x) T_t(x,y) \, \dtt, \qquad x,y \in X. }

Here we shall prove some preliminary estimate that will be needed later on.
\lem{a8}{
Suppose that \eqref{a3} -- \eqref{a6} are satisfied for $T_t$ and $\qq$. Then
\eqx{ \label{a8d} 
 \sup_{y\in X} \sum_{Q\in\qq} \int_{Q^{**}} \abs{R_j(x,y)} \abs{\psi_Q(x) - \psi_Q(y)} \, dx \leq C.
}
}

\begin{proof}
Fix $y\in X$ and $Q_0\in \qq$ such that  $y\in Q_0$. Write
\spx{
\sum_{Q\in\qq} \int_{Q^{**}} \abs{R_j(x,y)} & \abs{\psi_Q(x) - \psi_Q(y)} \, dx \leq  \sum_{Q\in\qq} \int_{Q^{**}} \abs{R^j_{d_{Q_0}, glob}(x,y)} \abs{\psi_Q(x) - \psi_Q(y)} \, dx \\
& + \sum_{Q\in\qq} \int_{Q^{**}\cap (Q_0^{***})^c} \abs{R^j_{d_{Q_0},loc}(x,y)} \abs{\psi_Q(x) - \psi_Q(y)} \, dx\\ 
& + \sum_{Q\in\qq} \int_{Q^{**}\cap Q_0^{***}} \abs{R^j_{d_{Q_0},loc}(x,y)} \abs{\psi_Q(x) - \psi_Q(y)} \, dx\\
& + \sum_{Q\in\qq} \int_{Q^{**}} \abs{R^j_V(x,y)} \abs{\psi_Q(x) - \psi_Q(y)} \, dx \\=  & S_1 + S_2 + S_3 + S_4.}

Using $\|\psi_Q\|_\8 \leq 1$, \eqref{finite_covering}, \eqref{a4},  \eqref{a3} and \eqref{a6} we have
\alx{
S_1 \sleq &  \int_X \int_{d_{Q_0}^2}^{\8} \abs{\partial_{x_j} T_t(x,y)} \, \dtt \, dx  \sleq 1,\\
S_{2} \sleq &   \int_{(Q_{0}^{***})^c} \int_{0}^{d_{Q_0}^2} \abs{\partial_{x_j} T_t(x,y)} \,\, \dtt \, dx \sleq 1,\\
S_4    \sleq &\int_X \int_0^\8 \abs{V_j(x)} T_t(x,y) \, \dtt \, dx \sleq 1.}

For $S_{3}$ consider $Q\in \qq$ such that $Q^{**}\cap Q_0^{***} \not= \emptyset$. The number of such $Q$ is bounded by an universal constant, $d_Q \simeq d_{Q_0}$, and $|\psi_Q(x) - \psi_Q(y)|\sleq d_{Q_0}^{-1} |x-y|$. Applying \eqref{a5} we obtain
\spx{
S_{3} \sleq &  \int_{Q_0^{***}} \int_{0}^{d_{Q_0}^2} \abs{\dxj \eee{  T_t(x,y) -   H_t(x-y)}} \, \frac{dt}{\sqrt{t}}  \, dx  \\
& +  \int_{Q_0^{***}}  \frac{|x-y|}{d_{Q_0}} \int_{0}^{d_{Q_0}^2} \abs{\dxj   H_t(x-y)} \, \frac{dt}{\sqrt{t}}  \, dx \\
&\sleq 1 +  \int_{Q_0^{***}}  \frac{|x-y|}{d_{Q_0}} \int_{0}^{\8} t^{-d/2} \exp\eee{-\frac{|x-y|^2}{ct}}\dttt \, dx\\
&\sleq 1 +  d_{Q_0}^{-1}  \int_{Q_0^{***}}  |x-y|^{-d+1} \, dx \sleq 1.
}
\end{proof}


\section{Proofs of Theorems \ref{thma} and \ref{thmb}.} \label{sec3}

\subsection{Proof of Theorem \ref{thma}} \label{proof-D}

\newcommand{\HarR}{H^1_{L, \rm{Riesz}}(X)}
\pr{
Denote
$$\norm{f}_{\HarR}:= \norm{f}_{L^1(X)} + \sum_{j=1}^d \norm{R_j f}_{L^1(X)}.$$

{\bf First inequality: $\norm{f}_{\HarR} \sleq \norm{f}_{\Har}$.} 
We shall show that
\eq{\label{rtyu}
\norm{R_j a}_{L^1(X)} \leq C
}
for $j=1,2,...,d$ and a $\qq$-atom $a$ with $C$ independent of $a$. In general, \eqref{rtyu} may not be enough to prove boundedness of an operator on $H^1$, see \cite{Bownik}. However, here Theorem~\ref{thm_decomp}, \eqref{rtyu}, and a standard continuity argument imply $\norm{f}_{\HarR} \sleq \norm{f}_{\Har}$. To show \eqref{rtyu}, according to Definition \ref{def_Qatoms}, suppose that $a$ is an $\qq$-atom associated with $Q\in \qq$. Let  $\RLQ$, $\RLQQ$ and $R^j_V$ denote the operators with the integral kernels defined in \eqref{loc_glob}. Applying \eqref{a6}, \eqref{a4}, \eqref{a3}, \eqref{a5},  and part~{\bf 1.} of Proposition \ref{prop_local_riesz} we have
\spx{ 
\norm{R_j a}_{L^1(X)} & \leq \norm{R^j_V a}_{L^1(X)} + \norm{\RLQQ a}_{L^1(X)} + \norm{\RLQ a}_{L^1((Q^{***})^c)} \\& + \norm{\left(\RLQ - \RDelta{j}{d_Q, loc}\right) a}_{L^1(Q^{***})}  + \norm{\RDelta{j}{d_Q,loc} a}_{L^1(Q^{***})} \leq C
}
and \eqref{rtyu} is proved. Let us notice here that since $a$ is bounded and $\supp \, a \subseteq Q^{**}$ then our assumptions guarantee that all the operators appearing above are well-defined.

{\bf Second inequality: $ \norm{f}_{\Har} \sleq \norm{f}_{\HarR}$.} Assume that $\norm{f}_{\HarR} <\8$. According to Theorem \ref{thm_decomp} it is enough to decompose $f$ as $\sum_k \la_k a_k$ with $\qq$-atoms $a_k$ and $\sum_k |\la_k| \leq \norm{f}_{\HarR}$. Let $\psi_Q$ be  a partition of unity related to $\qq$,  see Section \ref{sec_partition}. We have $f(x) = \sum_{Q \in \qq} f_Q(x)$, with $ f_Q = \psi_Q f$ and $\supp \, f_Q \subset \s{Q}{1}$. Notice that
\spx{
	\RDelta{j}{d_Q,loc} f_Q = & \left( \RDelta{j}{d_Q,loc} - \RLQ  \right)f_Q + \left(R_j  f_Q - \psi_Q R_j f\right)\\
	&- \RLQQ f_Q - R^j_V f_Q +  \psi_Q R_j f .
}
We use \eqref{a5}, Lemma \ref{a8}, \eqref{a4}, \eqref{a6} getting
\spx{
\sumQ \norm{\RDelta{j}{d_Q,loc} f_Q}_{L^1(Q^{**})} \leq & \sumQ \left[ \norm{\left( \RDelta{j}{d_Q,loc} - \RLQ \right) f_Q}_{L^1(Q^{**})} + 
\norm{R_j f_Q - \psi_Q R_j f }_{L^1(Q^{**})}\right.\\
& \left.+  \norm{ \RLQQ f_Q }_{L^1(Q^{**})}   +  \norm{ R^j_V f_Q }_{L^1(Q^{**})} + \norm{\psi_Q R_jf}_{L^1(Q^*)}\right]\\
\sleq & \sumQ \norm{f}_{L^1(Q^*)} +  \norm{f}_{L^1(X)} + \sumQ \norm{R_j f}_{L^1(Q^*)}\\
\sleq & \norm{f}_{\HarR}, }
for every $j = 1, ... ,d$. Now we use part {\bf 2.} of Proposition \ref{prop_local_riesz} for each $f_Q$, getting $\la_{Q,k}$, $a_{Q,k}$ such that
$$f_Q = \sum_k \la_{Q,k} a_{Q,k}, \qquad \sum_k |\la_{Q,k}| \sleq \norm{\RDelta{j}{d_Q,loc} f_Q}_{L^1(Q^{**})}. $$
The proof is finished by noticing that all $a_{Q,k}$ are $\qq$-atoms and
\spx{
f = \sum_{Q,k} \la_{Q,k} a_{Q,k}, \qquad \sum_{Q,k} \abs{\la_{Q,k}} \sleq \sumQ \norm{\RDelta{j}{d_Q,loc} f_Q}_{L^1(Q^{**})} \sleq \norm{f}_{\HarR} .
}
}

\subsection{Proof of Theorem \ref{thmb}} \label{proof-E} 
\pr{
According to Theorem \ref{thma} it is enough to prove \eqref{a0}--\eqref{a6} for the kernel
$$T_t(x,y) = \tti{1}(x_1,y_1)\cdot ... \cdot \tti{N}(x_N,y_N)$$
with the covering $\qq_1 \boxtimes ... \boxtimes \qq_N$, see \cite[Def.\ 1.5]{KPP}. It is enough to consider $N=2$ and then use an inductive argument. Assume that the conditions \eqref{a0}--\eqref{a6} are satisfied for $\tti{1}(x_1,y_1)$ and $\tti{2}(x_1,y_1)$ with $\qq_1$ and $\qq_2$, respectively.  The estimate \eqref{a0} for $T_t(x,y)$ follows directly. Moreover, \eqref{a1}--\eqref{a2} were already proved in the proof of  \cite[Thm.\ B]{KPP}.

To deal with \eqref{a3}--\eqref{a6} denote  $$\bx= (x_1,...,x_{d_1}, x_{d_1+1},...,x_{d_1+d_2})=(\bx_1,\bx_2) \in X_1\times X_2\subseteq \RR^{d_1}\times \RR^{d_2}.$$
Recall that a cuboid in $\qq_1 \boxtimes \qq_2$ is of the form $K=K_1\times K_2$, where $K_j\subseteq Q_j \in \qq_j$, $j=1,2$, and $d_K \simeq d_{K_1}\simeq d_{K_2} \simeq \min(d_{Q_1},d_{Q_2})$, see \cite[Def.\ 1.5]{KPP}. For the rest of the proof we fix $y\in K^{**} = K_1^{**}\times K_2^{**} \subseteq Q_1^{**}\times Q_2^{**}$ and without loss of generality we consider $\dxj $ for $j\in \set{d_1+1,...,d_1+d_2}$.

{\bf Proof of \eqref{a3}.} Notice that $(K^{***})^c = (K_1^{***}\times K_2^{***})^c = S_1 \cup S_2 \cup S_3 $, where
\eqx{
S_1 = X_1 \times (Q_2^{***})^c, \quad 
S_2 = X_1 \times (Q_2^{***} \setminus K_2^{***}), \quad 
S_3 = (K_1^{***})^c \times K_2^{***} .	
}
Using \eqref{a0} for $\tti{1}$ and \eqref{a3prim} for $\tti{2}$ we have
\spx{
\int_{S_1} \int_0^{d_K^2} \abs{\dxj  T_t(\bx,\by)} \frac{dt}{\sqrt{t}} \, d\bx & = \int_{S_1} \int_{0}^{d_K^{2}} \tti{1}(\bx_1,\by_1) \abs{\dxj  \tti{2}(\bx_2,\by_2) } \, \frac{dt}{\sqrt{t}} \, d\bx \\
& \leq \int_{(Q_2^{***})^c} \int_{0}^{c d_{Q_2}^{2}}  \abs{\dxj  \tti{2}(\bx_2,\by_2) } \, \frac{dt}{\sqrt{t}} \, d\bx_2 \sleq 1.
}

Using \eqref{a0} for $\tti{1}$ we have
\spx{
 \int_{S_2} \int_0^{d_K^2} \abs{\dxj  T_t(\bx,\by)} \frac{dt}{\sqrt{t}} \, d\bx \leq & \int_{S_2} \int_{0}^{d_K^{2}} \tti{1}(\bx_1,\by_1) \abs{\dxj  \tti{2}(\bx_2,\by_2) } \frac{dt}{\sqrt{t}}\, d\bx\\
 \leq &\int_{Q_2^{***}} \int_{0}^{d_K^{2}} \abs{\dxj \eee{  \tti{2}(\bx_2,\by_2) -  H_t(\bx_2-\by_2)}} \frac{dt}{\sqrt{t}} \, d\bx_2\\ 
 &+\int_{Q_2^{***}\setminus K_2^{***}} \int_{0}^{d_K^{2}} \abs{ \dxj  H_t(\bx_2-\by_2)} \frac{dt}{\sqrt{t}}\,  d\bx_2\\
 =&A_1+A_2.
 }
We have that $d_K \sleq d_{Q_2}$ and \eqref{a5prim} for $\tti{2}$ implies  $A_1\sleq 1$. Moreover, for $\by_2 \in K_2^{**}$ and $\bx_2\not\in K_2^{***}$ we have $|\bx_2-\by_2|\sgeq d_K$ and
\spx{
A_2 
& \sleq \int_{Q_2^{***} \setminus K_2^{***}} \int_{0}^{d_K^{2}} t^{-d_2/2} \exp\eee{-\frac{|\bx_2 - \by_2|^2}{ct}} \frac{dt}{t} \, d\bx_2 \\
&\sleq \int_{0}^{d_K^{2}} t^{M-d_2/2 -1} \, dt \cdot \int_{( K_2^{***})^c} \abs{\bx_2 - \by_2}^{-2M} \, d\bx_2 \sleq 1,
}
where $M$ is a fixed constant larger than $d_2/2$. What is left is to estimate the integral on $S_3$. Write
\spx{
 \int_{S_3} \int_0^{d_K^2} & \abs{\dxj  T_t(\bx,\by)} \frac{dt}{\sqrt{t}} \, d\bx 
\leq A_3 + A_4,}
where
\alx{
A_3 &= \int_{S_3} \int_{0}^{d_K^{2}} \tti{1}(\bx_1,\by_1) \abs{ \dxj \tti{2}(\bx_2,\by_2) - \dxj  H_t(\bx_2 -\by_2)} \frac{dt}{\sqrt{t}} \, d\bx,\\
A_4& = \int_{S_3} \int_{0}^{d_K^{2}} \tti{1}(\bx_1,\by_1) \abs{ \dxj  H_t(\bx_2 -\by_2)} \frac{dt}{\sqrt{t}} \, d\bx.
}
From \eqref{a0} for $\tti{1}$ and \eqref{a5prim} for $\tti{2}$ we easily get $A_3\sleq 1$. Let $\delta>0$ be fixed, Then,
\spx{
A_4 =& \int_{(K_1^{***})^c} \int_{K_2^{***} } \int_{0}^{d_K^{2}} t^{-2\de}\tti{1}(\bx_1,\by_1) \abs{ t^{\de+1/2} \dxj  H_t(\bx_2-\by_2)} \frac{dt}{t^{1-\delta}} \, d\bx_2 d\bx_1 \\
\sleq& \int_{(K_1^{***})^c} \sup_{s\leq d_K^2}\eee{ s^{-d_1/2 -2\de} \exp\eee{-\frac{|\bx_1 - \by_1|^2}{cs}}}\, d\bx_1 \\
&\times  \int_{K_2^{***} }\sup_{r\leq d_K^2}\eee{  r^{-d_2/2+\de} \exp\eee{-\frac{|\bx_2 - \by_2|^2}{cr}}}\, d\bx_2   \cdot \int_{0}^{d_K^{2}} t^{-1+\de} \, dt\\
\sleq &   \int_{|\bx_1-\by_1|\sgeq d_K} |\bx_1 - \by_1|^{-d_1-4\de}\, d\bx_1 \cdot \int_{|\bx_2-\by_2|\sleq d_K} |\bx_2 - \by_2|^{-d_2+2\de}\, d\bx_2 \cdot d_K^{2\de}\\
\sleq &  d_K^{-4\de} d_K^{2\de} d_K^{2\de} \sleq 1. 
}

{\bf  Proof of \eqref{a4}.} We have that $d_K \simeq d_{Q_1}$ or $d_K \simeq d_{Q_2}$. In the latter case $d_{K}\simeq d_{Q_2}$ the inequality \eqref{a4} for $T_t(x,y)$ follows simply from \eqref{a0} for $\tti{1}$ and \eqref{a4prim} for $\tti{2}$. Assume then that $d_K\simeq d_{Q_1} \sleq d_{Q_2}$. Let $t\geq d_K^2$ and $\by\in K^{**} \subseteq Q^{**}$. Write
\spx{
	\int_X \int_{d_K^2}^\8 \abs{\dxj  T_t(\bx,\by)} \frac{dt}{\sqrt{t}} \, d\bx = \int_X \int_{d_K^2}^{d_{Q_2}^2} .... + \int_X \int_{d_{Q_2}^2}^\8 ... = A_5 + A_6.
}

By \eqref{a0} for $\tti{1}$ and \eqref{a4} for $\tti{2}$ we easily get $A_6 \sleq 1$. Let $\de \in(0,\gamma)$ be as in \eqref{a1}--\eqref{a2}. For $A_5$ write
\spx{
A_5 \leq \int_{X_1} \sup_{t\geq d_K^2} \eee{ t^{\de} \tti{1}(\bx_1,\by_1)}\, d\bx_1 \cdot
 \int_{X_2} \int_{d_K^2}^{d_{Q_2}^2}  t^{-\de}\abs{\dxj  \tti{2}(\bx_2,\by_2)} \, \frac{dt}{\sqrt{t}}	\, d\bx_2  = A_{5,1} \cdot A_{5,2}.
	} 


By \eqref{a0} and \eqref{a1} for $\tti{1}$ we have
\spx{
A_{5,1} \sleq & \int_{Q_1^{***}} \sup_{t\geq d_K^2}  t^{\de-d_1/2} d\bx_1 + \int_{(Q_1^{***})^c} \sup_{t>0} t^{\de} \tti{1}(\bx_1,\by_1) \, d\bx_1\\
\sleq & d_{Q_1}^{d_1} d_K^{-d_1+2\de} +  d_{Q_1}^{2\de} \simeq d_K^{2\de }.
}
Moreover,
\spx{
A_{5,2} &\leq \int_{(Q_2^{***})^c} \int_{d_K^2}^{d_{Q_2}^2}  t^{-\de}\abs{\dxj  \tti{2}(\bx_2,\by_2)} \, \frac{dt}{\sqrt{t}}	\, d\bx_2 + \int_{Q_2^{***}} \int_{d_K^2}^{d_{Q_2}^2}  t^{-\de}\abs{\dxj  H_t(\bx_2-\by_2)} \, \frac{dt}{\sqrt{t}}	\, d\bx_2\\
&+\int_{Q_2^{***}} \int_{d_K^2}^{d_{Q_2}^2}  t^{-\de}\abs{\dxj  \tti{2}(\bx_2,\by_2)- \dxj  H_t(\bx_2-\by_2)} \, \frac{dt}{\sqrt{t}}	\, d\bx_2 = A_{5,2,1}+A_{5,2,2}+A_{5,2,3}.
}
Using \eqref{a3} and \eqref{a5} for $\tti{2}$ and the estimate $t^{-\de} \leq d_K^{-2\de}$ we easily get $A_{5,2,1}+A_{5,2,3} \sleq d_K^{-2\de}$. Also,
\spx{
A_{5,2,2} \leq \int_{d_K^2}^\8 t^{-1-\de} \int_{X_2} t^{-d_2/2} \exp\eee{-\frac{|\bx_2-\by_2|^2}{ct}}\, d\bx_2 \, dt \sleq  d_K^{-2\de}.
}
Combining all the estimates above we finish the proof of \eqref{a4} by noticing that $A_5+A_6 \sleq 1$.
	
{\bf Proof of \eqref{a5}}. We have that $d_K \simeq \min(d_{Q_1},d_{Q_2})$ and $K_j \subseteq Q_j$ for $j=1,2$. Using the triangle inequality write
\spx{
	&  \int_{K^{***}} \int_0^{d_K^2} \abs{\dxj\eee{  T_t(\bx,\by) -  H_t(\bx-\by)}} \, \frac{dt}{\sqrt{t}} \, d\bx  \\
	& \leq \int_{K^{***}} \int_0^{d_{K}^2}  \tti{1}(\bx_1,\by_1) \abs{\dxj\eee{  \tti{2}(\bx_2,\by_2) -   H_t(\bx_2-\by_2)}}  \, \frac{dt}{\sqrt{t}}\, d\bx \\
	& +  \int_{K^{***}} \int_0^{d_
	{K}^2} \abs{\tti{1}(\bx_1,\by_1) - H_t(\bx_1-\by_1)} \abs{\dxj  H_t(\bx_2-\by_2)}  \dtt \, d\bx  \\
	& = A_{7} + A_{8} \\
}

By \eqref{a0} for $\tti{1}$  and \eqref{a5prim} for $\tti{2}$ we have that $A_7\sleq 1$.

For $A_8$ we use \eqref{a2prim} for $\tti{1}$ obtaining
\spx{	A_{8}  \sleq &\int_{Q_1^{***}} \sup_{s\sleq d_{Q_1}^2} s^{-\de}\abs{\tsi{1}(\bx_1,\by_1) - H_s(\bx_1-\by_1)} \, d\bx_1  \\
	&\times \int_0^{cd_{Q_1}^2} \int_{Q_2^{***}} t^\de \abs{\dxj  H_t(\bx_2-\by_2)} \, d\bx_2 \, \frac{dt}{\sqrt{t}}\\ 
	\sleq & d_{Q_1}^{-2\de} \cdot \int_0^{d_{Q_1}^2} t^{-1+\de} \int_{X_2}   t^{-d_2/2} \exp\left( -\frac{|\bx_2-\by_2|^2}{c t} \right)  \, d\bx_2 \, dt \sleq 1.}

{\bf Proof of \eqref{a6}}. Fix $\by\in X$. Using \eqref{a0} for $\tti{1}$ and \eqref{a6} for $\tti{2}$ we have
\spx{ \int_{X} \int_0^\8 \abs{V_j(\bx)T_t(\bx, \by)} \, \dtt \, d\bx & \sleq \int_{X_2} \int_0^\8 \abs{V_j(\bx_2)}\tti{2}(\bx_2, \by_2) \int_{X_1} \tti{1}(\bx_1,\by_1) \, d\bx_1 \, \dtt \, d\bx_2 \\
& \sleq 1.
}
The proof of Theorem \ref{thmb} is finished.}

\section{Products of local and nonlocal atomic Hardy spaces} \label{sec3.1}

In this section we present an alternative version of Theorem \ref{thmb}. Consider the operator $L=-\Delta + L_2$, where $L_2$ is related to an admissible covering $\qq_2$ and satisfies all the assumptions of Section \ref{sec_assumptions}. It turns out that our methods work equally fine in this context. Notice that the Hardy space related to $-\Delta$ does not have local nature, but the Hardy space for $L=-\Delta + L_2$ will have local character as in Definition \ref{def_Qatoms}. Let us mention that this section will be needed to investigate the Dirichlet Laplacian on $\RR^d_+$, see Sections \ref{ssec_1.3} and \ref{sec:ex-dir-lap}.



More precisely, let $X = \RR^{d_1} \times X_2  $, where $X_2 \subseteq \RR^{d_2}$ is as in Section \ref{sec_assumptions}. We consider an operator $L_2$ on $L^2(X_2)$ and its semigroup $T_t^{[2]}$ with the kernel $T_t^{[2]}(x_2,y_2)$, $x_2,y_2 \in X_2$, $t>0$. Assume that $L_2$ and an admissible covering $\qq_2$ of $X_2$ satisfy \eqref{a0}--\eqref{a6}, see Section \ref{sec_assumptions} and Definition \ref{def_covering}. On $\RR^{d_1}$ we consider the Laplacian $-\Delta$ with the heat semigroup kernel $H_t(x_1-y_1)$, $x_1,y_1 \in \RR^{d_1}$, $t>0$, see \eqref{heat_kernel}. Following \cite[Sec. 1.4.4.]{KPP} we define the covering $\RR^{d_1} \boxtimes \qq_2  $ by splitting the strips $\RR^{d_1} \times Q_2$, $Q_2 \in \qq_2$, into countably many cuboids $Q_{1,n} \times Q_2$ such that $d_{Q_{1,n}}=d_{Q_2}$. Then $L=-\Delta+L_2$ is understood in the sense as in Section \ref{sec_products}.

The atomic characterization for $H^1_L(X)$ is given in  \cite[Cor. 1.14]{KPP}, where the atoms are related to the covering $\RR^{d_1}\boxtimes \qq_2$. In Theorem \ref{thm_riesz_loc_nonloc} below we provide a characterization by means of the Riesz transforms $R_j = D_j L^{-1/2}$, where  $D_j = \partial_{x_j}$  for $j=1,...,d_1$, and $D_j = \partial_{x_j} + V_j$ for $j=d_1+1,...,d_1+d_2$.
\thm{thm_riesz_loc_nonloc}{
Let $L=-\Delta + L_2$, where $-\Delta$ is the standard Laplacian on $\RR^{d_1}$ and $L_2$ with an admissible covering $\qq_2$ of $X_2\subseteq \RR^{d_2}$ satisfy \eqref{a0}--\eqref{a6}. Then there exists $C>0$ such that
\eqx{
C^{-1} \norm{f}_{\Har} \leq \norm{f}_{L^1(X)} + \sum_{j=1}^d \norm{R_j f}_{L^1(X)} \leq C  \norm{f}_{\Har}.
}
}
The proof of Theorem \ref{thm_riesz_loc_nonloc} follows directly from Theorem \ref{thma} and the following lemma.

\lem{lem_loc_nonloc}{
If \eqref{a0}--\eqref{a6} are satisfied for $L_2$ with an admissible covering $\qq_2$ of $X_2\subseteq \RR^{d_2}$, then \eqref{a0}--\eqref{a6} are satisfied for $L = -\Delta+L_2$ (see Section \ref{sec_products}) with an admissible covering $\RR^{d_1} \boxtimes \qq_2$ of $X = \RR^{d_1} \times X_2\subseteq \RR^{d_1+d_2}$.
}
Te proof of Lemma \ref{lem_loc_nonloc} uses the same techniques as the proof of Theorem \ref{thmb} and will be omitted.

\section{Examples} \label{sec4}
The goal of this section is to prove  Theorems \ref{coro_bess-riesz}, \ref{coro_lag}, and \ref{coro_dir}. According to Theorem \ref{thmb} it is enough to prove \eqref{a0}--\eqref{a6} for the one-dimensional Bessel operator $L_B^{[\be]}$ and the one-dimensional  Laguerre operator $\lla$.

Recall that \eqref{a0}--\eqref{a2} were proved in \cite[Prop.\ 4.3 and 4.5]{KPP}, so we shall deal only with \eqref{a3}--\eqref{a6} in Propositions \ref{riesz_bessel} and \ref{a_lag_riesz}. By $\Tt(x,y)$ we will denote the semigroup kernel related to: $L_B$ in Section \ref{sec:ex-bessel},  $L_L$ in Section \ref{sec:ex-laguerre}, and  $L_D$ in Section \ref{sec:ex-dir-lap}. Denote $\partial_x = \frac{d}{dx}$, the partial derivative on $(0,\8)$.

\subsection{Bessel operator.} \label{sec:ex-bessel}

The semigroup $T_{t}= \exp(-t \lbl)$ is given in terms of the integral kernel
\eq{\label{bessel-kernel}
T_{t}(x,y) = \frac{(xy)^{1/2}}{2t} I_{\be-1/2}\left(\frac{xy}{2t} \right) \exp\left(-\frac{x^2+y^2}{4t} \right), \qquad x,y\in X, t>0,
}
i.e. $T_{t} f(x) = \int_X T_{t}(x,y)f(y)\, dy$. Here, $I_\tau$ is the modified Bessel function of the first kind. For further reference recall some properties of the Bessel function $I_\tau$:
\begin{align} \label{Bessel-function-small} 
&&I_\tau(x) & = C_\tau x^\tau + O(x^{\tau+1}), & \text{ for } & x\sim 0,&&\\
\label{Bessel-function-large}
&&I_\tau(x) & = (2\pi x)^{-1/2} e^x + O(x^{-3/2}e^x), & \text{ for } & x \sim \8,&&\\
\label{der-bessel}
&&\partial_x(x^{-\tau}I_\tau(x)) &= x^{-\tau} I_{\tau+1}(x)& \text{ for } & x >0,&&
\end{align}
see e.g. \cite{Watson}. The main goal of this section is to prove the following proposition.
\prop{riesz_bessel}{
Let $X=(0,\8)$ and $\be>0$. Then \eqref{a3}--\,\eqref{a6} hold for $L_B^{[\be]}$ with $\qq_B$, see \eqref{Bess_cov}.
}

\pr{Using \eqref{der-bessel} we have
\sp{ \label{dxBessel}
\partial_x T_t(x,y) = \frac{(xy)^{\frac{1}{2}}}{2t} \exp\left(-\frac{x^2+y^2}{4t} \right) \left( \frac{y}{2t} I_{\beta+\frac{1}{2}}\left(\frac{xy}{2t} \right) + \frac{\beta}{x} I_{\beta-\frac{1}{2}}\left(\frac{xy}{2t} \right) - \frac{x}{2t} I_{\beta-\frac{1}{2}}\left(\frac{xy}{2t} \right)  \right).
}
Denote {\bf case 1}: $xy \sleq t$. In this case, by \eqref{dxBessel} and  \eqref{Bessel-function-small},
	\sp{\label{diff1}
	\abs{\partial_x T_t(x,y)} & \sleq t^{-1/2} \eee{\frac{xy}{t}}^{\be}  \exp\left( -\frac{x^2+y^2}{ct} \right) \eee{\frac{1}{x}+\frac{x}{t}}.
	}
In {\bf case 2}: $t \sleq xy$, using \eqref{dxBessel} and  \eqref{Bessel-function-large}, we have
	\sp{\label{diff2}
	\abs{\partial_x T_t(x,y)} & \sleq \frac{x+y}{t^{3/2}} \exp\left( -\frac{|x-y|^2}{ct} \right).
	}
For the rest of the proof let us fix $I = [2^n, 2^{n+1}] \in \qq_B$ and $y\in I^{**}$. Then $y \simeq 2^n = d_{I}$. Fix $2^{-1}<\kappa_1<1<\kappa_2<2$ such that $I^{***} = [\kappa_1 2^n, \kappa_2 2^{n+1}]$.

{\bf Proof of \eqref{a3}.} Write 
\spx{\int_{(I^{***})^c} \int_0^{d_I^2} \abs{\partial_x T_t(x,y)} \frac{dt}{\sqrt{t}} \, dx  \leq & \int_{0}^{\kappa_1 2^n} \int_0^{xy} ... + \int_{0}^{\kappa_1 2^n} \int_{xy}^{2^{2n}} ... + \int_{\kappa_2 2^{n+1}}^{\8} \int_0^{2^{2n}} ...
 \\
 = & A_1 + A_2 + A_3.
 }
 
 
For $A_1$ and $A_3$ we use \eqref{diff2}, whereas for $A_2$ we use \eqref{diff1}, obtaining:
\alx{
A_1 & \sleq \int_0^{\kappa_1 2^{n}} \int_0^{xy} \frac{2^{n}}{t^{3/2}} \exp\eee{-\frac{2^{2n}}{ct}} \dtt \, dx\sleq 2^{-n} \int_0^{2^n} \int_0^{2} \exp\eee{-\frac{1}{ct^2}} \, \frac{dt}{t^2}\, dx \sleq 1,\\
A_2 &\sleq  \int_{0}^{\kappa_1 2^{n}} \int_{xy}^{2^{2n}} \eee{\frac{x2^{n}}{t}}^\be  \exp\eee{-\frac{2^{2n}}{ct}} \eee{\frac{1}{x}+\frac{x}{t}} \frac{dt}{t}\, dx\\
    &\sleq  2^{- n\be} \cdot \int_{0}^{2^{n}} x^{-1+\be}\, dx \cdot  \int_{0}^{\8} \eee{\frac{2^{2n}}{t}}^\be \exp\eee{-\frac{2^{2n}}{ct}}  \eee{1+\frac{2^{2n}}{t}} \frac{dt}{t} \sleq 1,\\
A_3 & \sleq \int_{\kappa_2 2^{n+1}}^\8 \int_0^{2^{2n}} \frac{x}{t^{3/2}} \exp\eee{-\frac{x^2}{ct}} \dtt \, dx \\
&\sleq  \int_{2^{n+1}}^\8 x^{1-2N}\, dx \cdot \int_0^{2^{2n}} t^{N-2}\, dt \sleq  1,
}
where $N$ is arbitrarily large constant (here $N>1$ is enough).

{\bf Proof of \eqref{a4}.} Let us write
\spx{\int_X \int_{d_I^2}^\8 \abs{\partial_{x} T_t(x,y)} \frac{dt}{\sqrt{t}} \, dx  & =  \int_0^{2^{n+2}} \int_{2^{2n}}^\8 ... + \int_{2^{n+2}}^\8 \int_{2^{2n}}^{2^nx} ... + \int_{2^{n+2}}^\8 \int_{2^n x}^\8 ... \\
&= A_4+A_5+A_6.
}
For $A_4$ we have observe that  $x/t \sleq  2^n/t \sleq x^{-1}$. Using \eqref{diff1},
\spx{
A_4 \sleq &  \int_{0}^{2^{n+2}} \int_{2^{2n}}^\8 \eee{\frac{x2^{n}}{t}}^\be  \exp\eee{-\frac{2^{2n}}{ct}} \eee{\frac{1}{x}+\frac{x}{t}} \frac{dt}{t}\, dx\\
\sleq &  2^{n\be}\int_{0}^{2^{n+2}}  x^{-1+\be} dx \cdot \int_{2^{2n}}^\8 t^{-1-\be}  \exp\eee{-\frac{2^{2n}}{ct}} \, dt \sleq 1.\\
}
In $A_5$ and $A_6$ we use \eqref{diff2} and \eqref{diff1}, respectively. For an arbitrary large $N$ we have:
\alx{
A_5 \sleq & \int_{2^{n+2}}^\8 \int_{2^{2n}}^{2^n x} \frac{x}{t^{3/2}}  \exp\eee{-\frac{x^2}{ct}} \dtt\, dx\\
\sleq &  \int_{2^{n+2}}^\8  x^{-2N+1} \int_{0}^{2^n x} t^{N-2}  \, dt \, dx \sleq 1,\\
A_6 \sleq &  \int_{2^{n+2}}^\8 \int_{2^{n}x}^{\8} \eee{\frac{2^n x}{t}}^{\be}  \exp\eee{-\frac{x^2}{ct}} \frac{1}{x} \eee{1+\frac{x^2}{t}}\,\frac{dt}{t}\, dx\\
\sleq &  2^{n\be} \int_{2^{n+2}}^\8  x^{-1-\be} \int_{0}^{\8} \eee{\frac{x^2}{t}}^{\be} \eee{1+\frac{x^2}{t}} \exp\eee{-\frac{x^2}{ct}}  \, \frac{dt}{t}\, dx\\
\sleq &  2^{n\be} \int_{2^{n+2}}^\8  x^{-1-\be} \, dx \cdot \int_0^\8 t^{\be-1} (1+t) e^{-t}\, dt  \sleq 1.\\
}

{\bf Proof of \eqref{a5}.} Observe that for $x\in I^{***}$, $y\in I^{**}$, and $t\leq d_I^2 = 2^{2n}$ we have $t\sleq  xy$. Therefore, using \eqref{dxBessel} and \eqref{Bessel-function-large} we get
\sp{
\label{otherdxB} \partial_x T_t(x,y) & = \frac{y-x}{2t} \frac{(xy)^{1/2}}{2t} \exp\left(-\frac{|x-y|^2}{4t} \right) \left( \frac{\pi xy}{ t} \right)^{-1/2} + R(x,y) \\
	& = \partial_x H_t(x,y) + R(x,y),
	}
where 
\sp{\label{eq:410}
\abs{R(x,y)} &\sleq  t^{-1/2}  \exp\left( -\frac{|x-y|^2}{4t}  \right) \left( \frac{x+y}{xy} + x^{-1} \right) \\
	&\sleq  x^{-1} t^{-1/2}  \exp\left( -\frac{|x-y|^2}{4t}  \right),
}
since $x\simeq y \simeq d_I$. Notice that $|x-y|\sleq 2^n$. By \eqref{otherdxB} and \eqref{eq:410} we obtain
\spx{\int_{I^{***}} \int_0^{2^{2n}} & \abs{\partial_x T_t(x,y) - \partial_x H_t(x,y)} \frac{dt}{\sqrt{t}} \, dx \leq \int_{I^{***}} \int_0^{2^{2n}} \abs{R(x,y)} \frac{dt}{\sqrt{t}} \, dx \\
	&\sleq  \int_{I^{***}} x^{-1} \int_0^{2^{2n}} \exp\left( -\frac{|x-y|^2}{4t} \right) \frac{dt}{t} \, dx \\
	&= C \int_{I^{***}} x^{-1} \int_{|x-y|^2/2^{2n}}^\8 \exp\left( -t/4 \right) \frac{dt}{t} \, dx  \\
 	&\sleq  \int_{2^{n-1}}^{2^{n+2}} \ln \eee{2+\frac{2^n}{|x-y|}}  \, \frac{dx}{x} \sleq  \int_{-2}^2 \ln \eee{2+|x|^{-1}}\, dx \sleq 1.
}

{\bf Proof of \eqref{a6}.} Using \eqref{bessel-kernel}, \eqref{Bessel-function-small} and  \eqref{Bessel-function-large},  we have that
\spx{
\int_0^\8 T_t(x,y) \, \dtt   &\sleq
\int_0^{xy} \exp\eee{-\frac{|x-y|^2}{4t}} \frac{dt}{t} +  \int_{xy}^\8 \eee{\frac{xy}{t}}^{\beta} \exp\eee{-\frac{x^2+y^2}{ct}}\, \frac{dt}{t}\\
&\sleq \begin{cases} \eee{x/y}^\beta & \quad x \leq y/2, \\
 \ln\eee{y|x-y|^{-1}} & \quad |x-y|\leq y/2,\\
\eee{y/x}^\beta &  \quad x\geq 3y/2.
\end{cases}}
Hence,
\spx{
\int_X \int_0^\8 x^{-1} T_t(x,y) \, \dtt & \sleq y^{-\beta} \int_0^{y/2} x^{-1+\beta} \, dx + \int_{|x-y|\leq y/2} \ln\eee{\frac{y}{|x-y|}} \, \frac{dx}{x} \\
& + y^{\beta} \int_{3y/2}^\8 x^{-1-\beta} \, dx \sleq 1.
}
This ends the proof of Proposition \ref{riesz_bessel}.
}

\newcommand{\sh}{\rm{sh}}
\newcommand{\ch}{\rm{ch}}
\subsection{Laguerre operator.} \label{sec:ex-laguerre}

Recall that $\be > 0$ denotes the parameter related to the Lagurre operator $L_{L}^{[\be]}$, see \eqref{L_Lagg}. The goal of this section is to prove we have the following proposition.
\prop{a_lag_riesz}{
Let $X=(0,\8)$ and $\be>0$. Then \eqref{a3}--\,\eqref{a6} hold for $L_{L}^{[\be]}$ with $\qq_{L}$ given in \eqref{Lag_cov}.
}

Before going to the proof let us make some preparations. In what follows we shall use the notation $\sh(t) = \sinh(t)$, and $\ch(t) = \cosh(t)$. The semigroup $T_t = T_{L,t}=\exp\left(-tL_{L}^{[\be]}\right)$ has a kernel given by
\eq{\label{lag1-kernel}
T_{t}(x,y) = \frac{(xy)^{1/2}}{\sh (2t)} I_{\be-1/2}\left(\frac{xy}{\sh (2t)} \right) \exp\left(-\frac{\ch (2t)}{2\,{\sh{(2t)}}}(x^2+y^2) \right), \qquad x,y\in X, \ t>0.
}
Denote 
\eq{\label{U}
U_{\be-1/2}(x) = I_{\be-1/2}(x) \exp(-x)\sqrt{2\pi x},
}
so that
\eq{\label{U_prop}
 |U_{\be-1/2}(x) - 1|\sleq x^{-1}, \quad |U_{\be-1/2}(x) - U_{\be+1/2}(x)|\sleq x^{-1}, \qquad x\sim \8,
}
c.f. \eqref{Bessel-function-large}. Denote
\eqx{\Theta(t,x,y) = \exp\left( \frac{(1-\ch(2t))(x^2+y^2)}{2\sh(2t)}\right).}

In some cases we shall use different expression for $T_t(x,y)$, namely
\eq{\label{lag2-kernel}
T_{t}(x,y) = \frac{\Theta(t,x,y)}{\sqrt{2\pi \sh (2t)}} U_{\be-1/2}\left(\frac{xy}{\sh (2t)} \right) \exp\left(-\frac{|x-y|^2}{2\sh{(2t)}} \right), \qquad x,y\in X, \ t>0.
}
Using \eqref{lag1-kernel}, \eqref{lag2-kernel}, \eqref{der-bessel}, and \eqref{U} we get three expressions for $\partial_x T_t(x,y)$, i.e.
\al{ 
\partial_x T_t(x&,y) =
 \frac{\sqrt{xy}}{\sh(2t)} \exp\left( -\frac{\ch(2t)}{2\sh(2t)} \left(x^2+y^2\right) \right) \cdot F_1(t,x,y)
 \label{R1}
 \\
\label{R2}
 &= \frac{ \Theta(t,x,y)}{\sqrt{2\pi\sh(2t)}} \exp\left( -\frac{|x-y|^2}{2\sh(2t)} \right)  \cdot F_2(t,x,y)\\
\label{R3}
&= \frac{ \Theta(t,x,y)}{\sqrt{2\pi\sh(2t)}} \exp\eee{ -\frac{|x-y|^2}{2\sh(2t)} } \cdot\eee{ \frac{y-x}{\sh(2t)} U_{\be+1/2} \eee{\frac{xy}{\sh(2t)}}+F_3(t,x,y)}, }
where
\alx{
F_1(t,x,y) =  & \frac{y}{\sh(2t)}I_{\be+1/2}\left( \frac{xy}{\sh(2t)}\right) + \frac{\be}{x} I_{\be-1/2} \left( \frac{xy}{\sh(2t)}\right) - x \frac{\ch(2t)}{\sh(2t)} I_{\be-1/2} \left( \frac{xy}{\sh(2t)}\right),\\
  F_2(t,x,y) = & \frac{y}{\sh(2t)}U_{\be+1/2}\left( \frac{xy}{\sh(2t)
  }\right) + \frac{\be}{x} U_{\be-1/2} \left( \frac{xy}{\sh(2t)}\right) -x \frac{\ch(2t)}{\sh(2t)} U_{\be-1/2} \left( \frac{xy}{\sh(2t)} \right),\\
 F_3(t,x,y) =  & \frac{\be}{x} U_{\be-1/2} \left( \frac{xy}{\sh(2t)}\right) -\frac{x}{\sh(2t)}\left(\ch(2t) U_{\be-1/2} \left( \frac{xy}{\sh(2t)}\right) -  U_{\be+1/2} \left( \frac{xy}{\sh(2t)}\right) \right).
}

Observe that
\al{\label{Th1}
0 < & \Theta(t,x,y)\sleq  \exp\left( - c t(x^2+y^2) \right),  \qquad && \text{ for } \ t \sleq 1, \ x,y\in X
\\
\label{Th2}
0 < & \Theta(t,x,y)\sleq  \exp(-c(x^2 +y^2)), \qquad && \text{ for } \ t\sgeq 1, \  x,y\in X.
}
Moreover, using \eqref{Bessel-function-small} and \eqref{U_prop} we get
\al{ \label{F1}
|F_1(t,x,y)| \sleq &  \eee{\frac{xy}{\sh(2t)}}^{\be-1/2} \eee{ \frac{1}{x}+\frac{x\ch(2t)}{\sh(2t)}}, && xy\sleq \sh(2t),\\
\label{F2}
|F_2(t,x,y)| \sleq &  \eee{\frac{y}{\sh(2t)}+\frac{x \ch(2t)}{\sh(2t)}}, && xy\sgeq \sh(2t),\\
\label{F3}
|F_3(t,x,y)| \sleq &  \eee{\frac{1}{x}+xt + \frac{1}{y}}, && xy\sgeq \sh(2t), \  t\leq 1.
}

Now we are almost ready to prove Proposition \ref{a_lag_riesz} but first let us make a few comments and fix some notion. The proof relies on a detailed and lengthy analysis, but essentially one uses only simple calculus and properties of $I_{\be-1/2}$. We shall write $a\land b = \min(a,b)$ and $a \lor b = \max(a,b)$. Recall that $\qq_L$ is the set of intervals given in \eqref{Lag_cov}.
The proof will be given in two cases. First we shall deal with the  sub-intervals of $[0,1]$ in Section \ref{sssec1}. Then we shall consider sub-intervals of $[1,\8)$ in Section \ref{sssec2}. The letter $n$ will always be a positive integer. Moreover, we shall use $N$ as a constant that is fixed and large enough, depending on the context (most often we shall use the inequality $\exp(-x)\sleq  x^{-N}$).

\subsubsection{Case 1: $I\subseteq [0,1]$. }\label{sssec1} We consider $I=[2^{-n},2^{-n+1}]$, $n\in \NN$, and $y\in I^{**}$. Then $y\simeq 2^{-n} =d_I$. Fix $2^{-1}<\kappa_1<1<\kappa_2<2$ such that $I^{***} = [\kappa_1 2^{-n}, \kappa_2 2^{-n+1}]$.

{\bf Proof of \eqref{a3}  in Case 1.} We deal with $0< t \leq 2^{-2n} \leq 1$, $\sh(t) \simeq t$ and $\ch(t) \simeq 1$. Then
\spx{\int_{(I^{***})^c} \int_0^{d_I^2} \abs{\partial_x T_t(x,y)} \frac{dt}{\sqrt{t}} \, dx  \leq &\int_0^{\kappa_1 2^{-n}} \int_0^{2^{-2n}\land xy} ... + \int_{\kappa_2 2^{-n+1}}^{\8} \int_0^{2^{-2n}\land xy} ...\\&+ \int_{(I^{***})^c} \int_{2^{-2n} \land xy}^{2^{-2n}} ...
 = A_1 + A_2 + A_3.}



For $A_1$ we have $xy \sgeq t$, $x< y$, $|x-y|\simeq y$, and $|F_2(t,x,y)| \sleq y/t$. Using \eqref{R2}, \eqref{Th1}, and \eqref{F2},
\spx{A_1 
 \sleq & \  y \int_0^{2^{-n}}  \int_0^{2^{-2n}} t^{-1} \exp\left( -\frac{y^2}{ct} \right) \frac{dt}{t} \, dx \sleq y^{1-2N} \int_0^{2^{-n}} dx \cdot \int_0^{2^{-2n}} t^{N-2} \, dt \sleq 1. }
 
For $A_2$ we have $ xy \sgeq t$, $y< x$, $|x-y|\simeq x$, and $|F_2(t,x,y)| \sleq x/t$. Using \eqref{R2}, \eqref{Th1}, and \eqref{F2},
\spx{A_2 
 \sleq &  \int_{2^{-n+1}}^\8  x \int_0^{2^{-2n}} t^{-1} \exp\left( -\frac{x^2}{ct} \right) \frac{dt}{t} \, dx \sleq  \int_{2^{-n+1}}^\8 x^{1-2N} dx \cdot \int_0^{2^{-2n}} t^{N-2} \, dt \sleq 1. }

Notice that $A_3$ appears only when $x \leq \kappa_1 2^{-n}$. Moreover, $x^2 \sleq xy \sleq t$, and $|F_1(t,x,y)| \sleq  x^{-1} (xy/t)^{\be-1/2} $. Using \eqref{R1} and \eqref{F1},

\spx{A_3 & \sleq \int_0^{\kappa_1 2^{-n}} x^{-1} \int_{0}^{2^{-2n}} \left( \frac{xy}{t}\right)^{\be}  \exp\left( -\frac{y^2}{ct} \right) \frac{dt}{t} \, dx\\
& \sleq y^{-2N+\be} \int_0^{2^{-n}} x^{\be-1} \, dx \cdot \int_{0}^{2^{-2n}}t^{N-\be-1} \, dt \sleq 1.}

{\bf Proof of \eqref{a4}  in Case 1.} Recall that $y\simeq 2^{-n}$. We shall consider $t\geq d_I^2 = 2^{-2n}$. Write
\spx{\int_0^\8 \int_{d_I^2}^\8 \abs{\partial_x T_t(x,y)}\, \dtt \, dx =& \int_0^{2^{-n+3}} \int_{2^{-2n}}^1 ... + \int_{2^{-n+3}}^\8 \int_{2^{-2n}}^{1 \land xy} ... + \int_{2^{-n+3}}^\8 \int_{1 \land xy}^1 ... \\
& +  \int_{0}^\8 \int_1^{1 \vee \ln(\sqrt{xy})} ... + \int_{0}^\8 \int_{1\vee \ln(\sqrt{xy})}^\8 ... \\
=& A_4 + A_5 + A_6 + A_7 + A_8.
}
In the integrals $A_4$--$A_6$ we have $t\leq 1$, so that  $\sh(2t)\simeq t$ and $\ch(2t) \simeq 1$.

For $A_4$ we have $ x^2 \sleq xy \sleq t$, so that  $|F_1(t,x,y)| \sleq  x^{-1} (xy/t)^{\be-1/2} $. Using \eqref{R1} and \eqref{F1},
\eqx{
A_4 \sleq \int_0^{2^{-n+3}} x^{-1} \int_{2^{-2n}}^1 \eee{\frac{xy}{t}}^{\be} \, \dttt  \, dx \sleq y^{\be} \int_0^{2^{-n+3}} x^{\be-1} \, dx \cdot \int_{2^{-2n}}^\8 t^{-\be-1}\, dt \sleq 1.
}

For $A_5$ we have $ xy \sgeq t$ and $|x-y| \simeq x \geq y$ , since $x\geq 2^{-n+3}$ and $y\leq 2^{-n+2}$. Then $|F_2(t,x,y)| \sleq x/t$. Using \eqref{R2},  \eqref{Th1}, and \eqref{F2},
\spx{
A_5 &\sleq \int_{2^{-n+3}}^\8 x \int_{2^{-2n}}^{2^{-n+2}x} \exp\eee{-\frac{x^2}{ct}}  \, \dtttt  \, dx \sleq\int_{2^{-n+3}}^\8 x^{1-2N} \int_{0}^{2^{-n+2}x} t^{N-2}\, dt \, dx \sleq 1 .
}

For $A_6$ we have $ xy \sleq t$ and $ x \geq y$. Then $|F_1(t,x,y)| \sleq  x^{-1} (xy/t)^{\be-1/2}(1+x^2/t) $. Using \eqref{R1} and \eqref{F1},
\spx{
A_6 &\sleq \int_{2^{-n+3}}^\8 x^{-1} \int_{xy}^1 \eee{\frac{xy}{t}}^{\be}  \exp\eee{-\frac{x^2}{ct}}  \eee{ 1+\frac{x^2}{t}}  \, \dttt  \, dx \\
&\sleq y^{\be}\int_{2^{-n+3}}^\8 x^{\be-1 } \int_0^\8 t^{-\be-1} \exp\eee{-\frac{x^2}{c't}} \, dt \, dx \\
& \sleq y^{\be} \int_{2^{-n+3}}^\8 x^{-\be-1} \, dx \cdot \int_0^\8 t^{-\be-1}  \exp\eee{-\frac{1}{c't}} \, dt \sleq 1.
}

In the integrals $A_7$--$A_8$ we deal with $t>1$, so that $\sh(2t) \simeq e^{2t}$ and $\sh(2t)/\ch(2t) \simeq 1$.

The term $A_7$ appears only when $x \sgeq 2^n$. Here $xy\sgeq \sh(2t)$,  $x>y$, and $|F_2(t,x,y)| \sleq x $. Using \eqref{R2}, \eqref{Th2}, and \eqref{F2},
\spx{ 
A_7 &\sleq \int_{0}^\8 \int_1^{\8} \frac{ x}{(\sh(2t))^{1/2}} \exp\eee{- c x^2} \, \dtt \, dx \sleq 1.
}

For $A_8$ we have $ xy \sleq \sh(2t)$ and $|F_1(t,x,y)|\sleq (xy/\sh(2t))^{\be-1/2} (x+x^{-1})$. Using \eqref{R1} and \eqref{F1},
\spx{
A_8 &\sleq \int_0^{\8} \int_1^\8 \frac{(xy)^{\be}}{(\sh(2t))^{\be+1/2}}  \exp\eee{- cx^2} \eee{x+x^{-1}} \, \dtt \, dx \\
&\sleq y^{\be} \int_0^\8 x^{\be} \eee{x+x^{-1}} \exp\eee{- cx^2} \, dx \cdot \int_1^\8 (\sh(2t))^{-\be-1/2} \, \dtt \sleq 1,
}
where we have used that $y\leq 2$ and $\be>0$.

{\bf Proof of \eqref{a5}  in Case 1.}
In \eqref{a5} we deal with $x\simeq y \simeq 2^{-n}$ and $t\leq 2^{-2n}$, so $t \sleq xy \sleq 1$. Recall that $H_t(x-y)$ denotes the classical heat kernel on $\RR$. Using  \eqref{R3},
\sp{\label{k123}
&\abs{\partial_x (T_t(x,y) -  H_t(x-y))}
\\
\leq & \abs{\left(\partial_x H_{\frac{\sh(2t)}{2}}(x-y) - \partial_x H_t(x-y) \right)\Theta(t,x,y)U_{\be+\frac{1}{2}}\left( \frac{xy}{\sh(2t)}\right)}\\
& + \abs{\Theta(t,x,y)-1 - \Theta(t,x,y) \eee{1-U_{\be+1/2}\eee{\frac{xy}{\sh(2t)}}}} \cdot \abs{ \partial_x H_t(x,y)} \\
& + \frac{\Theta(t,x,y)}{(2\pi \sh(2t))^{1/2}} \exp\eee{-\frac{|x-y|^2}{2\sh(2t)}} F_3(t,x,y) \\
 = & K^{[1]}_{t}(x,y) + K^{[2]}_{t}(x,y) + K^{[3]}_{t}(x,y) .
}

Recall that $t\leq 1$ and notice that $
\abs{\partial_t \partial_x H_t(x-y)} \sleq t^{-3/2} \exp\eee{-|x-y|^2/(8t)}
$. Using \eqref{U_prop}, \eqref{Th1},  and the mean-value theorem we have
\spx{
K^{[1]}_t(x,y) & \sleq \abs{\sh(2t)/2 - t} t^{-3/2} \exp\left(-|x-y|^2/(ct)\right) \sleq t^{3/2}.
}
Therefore, 
\spx{ \label{k3}
 \int_{I^{***}} \int_0^{d_I^2} K_t^{[1]}(x,y) \frac{dt}{\sqrt{t}} \, dx & \sleq \int_{2^{-n-1}}^{2^{-n+2}} \, dx \cdot  \int_0^{2^{-2n}}  t\, dt \sleq 1.
}

Turning to $K_t^{[2]}$ notice that
\sp{ \label{theta-1}
\abs{1 -  \Theta(t,x,y) }  & = \abs{ \exp(0) -\exp\left( \frac{(1-\ch(2t)) (x^2+y^2)}{2\sh(2t)} \right)} \sleq 
t y^2.
}
Using \eqref{U_prop}, \eqref{Th1} and \eqref{theta-1} we get $K^{[2]}_t(x,y) \sleq t (y^2+(xy)^{-1}) |\partial_x H_t(x-y)| \sleq  2^{2n} $, hence
\spx{
 \int_{I^{***}} \int_0^{2^{-2n}} K_t^{[2]}(x,y) \frac{dt}{\sqrt{t}} \, dx & \sleq 2^{2n} \cdot \int_{2^{-n-1}}^{2^{-n+2}}dx \cdot  \int_0^{2^{-2n}}\, \dtt \sleq 1.
 }

For $K^{[3]}_t$ by \eqref{F3} we have $|F_3(t,x,y)| \sleq y^{-1} \sleq 2^n$. Using \eqref{Th1},
\spx{
 \int_{I^{***}} \int_0^{d_I^2} K_t^{[3]}(x,y) \frac{dt}{\sqrt{t}} \, dx & \sleq 2^n \int_{2^{-n-1}}^{2^{-n+2}}  \int_{0}^{2^{-2n}} \exp\left( -\frac{|x-y|^2}{ct}\right) \frac{dt}{t} \, dx  \\
 & \sleq 2^n \int_{|x-y|\sleq 2^{-n}} \int_{2^{2n} |x-y|^2}^{\infty} e^{-t}\, \dttt \, dx \\
 & \sleq 2^n \int_{|x-y|\sleq 2^{-n}} \ln (2^{-n} |x-y|^{-1}) \, dx \sleq \int_{|x|\sleq 1 } \ln |x|^{-1}\, dx \sleq  1.
}

\subsubsection{Case 2: $I\subseteq [1,\8)$. } \label{sssec2} Fix $y\in I^{**}$ and $n\in \NN$ such that $I\subseteq[2^{n},2^{n+1}]$. We have $y\simeq 2^{n} =d_I^{-1}$.

{\bf Proofof \eqref{a3} in Case 2.}

Notice that we deal with $0< t \leq 2^{-2n} \leq 1$, $\sh(t) \simeq t$ and $\ch(t) \simeq 1$. For $y\in I^{**}$ and $x\not\in I^{***}$ we have $|x-y|\sgeq 2^{-n}$, so that
\spx{\int_{(I^{***})^c} \int_0^{2^{-2n}} \abs{\partial_x T_t(x,y)} \frac{dt}{\sqrt{t}} \, dx  \leq &\int_{0}^\8 \int_{2^{-2n} \land xy}^{2^{-2n}} ... + \int_{2^{-n} \sleq |x-y| \leq 2^{n-2}} \int_0^{2^{-2n}\land xy} ...\\
&+ \int_{ |x-y| \geq 2^{n-2}} \int_0^{2^{-2n}\land xy} ...
 = A_9 + A_{10} + A_{11}.}

For $A_9$ we have $xy \sleq t$ and $x\sleq 2^{-3n}$, so that $|x-y| \simeq y$. Thus $|F_1(t,x,y)| \sleq x^{-1}(xy/t)^{\be-1/2}$. Using \eqref{R1} and \eqref{F1},
\spx{A_9 & \sleq \int_0^{c2^{-3n}}  x^{-1} \int_{0}^{2^{-2n}} \left( \frac{xy}{t}\right)^{\be}  \exp\left( -\frac{y^2}{ct} \right) \dttt \, dx\\
& \sleq y^{-2N+\be} \int_0^{c 2^{-3n}} x^{\be-1} \, dx \cdot \int_{0}^{2^{-2n}}t^{N-\be-1} \, dt \sleq 2^{-4Nn} \sleq  1.}

For $A_{10}$ we have $ xy \sgeq t$, $x \simeq y \simeq 2^n$, $x^{-1} \sgeq xt$, so that $|F_3(t,x,y)| \sleq y^{-1} \sleq |x-y|/t$. Using \eqref{R3}, \eqref{F3}, and \eqref{Th1},
\spx{A_{10} 
 \sleq & \int_{2^{-n} \sleq |x-y| \leq 2^{n-2}} |x-y| \int_0^{2^{-2n}} \exp\left( -\frac{|x-y|^2}{ct} \right) \dtttt \, dx\\
 \sleq &  \int_{2^{-n} \sleq |x-y|} |x-y|^{1-2N} \, dx \cdot \int_0^{2^{-2n}} t^{N-2}\, dt  \sleq 1. }

For $A_{11}$ we have $ xy \sgeq t$, $|x-y| \simeq x+y$, and $|F_2(t,x,y)| \sleq (x+y)/t$. Using \eqref{R2}, \eqref{Th1}, and \eqref{F2},
\spx{A_{11} & \sleq \int_{|x-y|\geq 2^{n-2}} (x+y) \int_0^{2^{-2n}} \exp\left( -\frac{(x+y)^2}{ct} \right) \dtttt \, dx\\
& \sleq \int_{0}^\8 (x+y)^{1-2N}  \, dx \cdot \int_0^{2^{-2n}} t^{N-2}\, dt \sleq 2^{4n(1-N)} \sleq  1.}

{\bf Proof of \eqref{a4} in Case 2.}

Write
\spx{\int_{X} \int_{d_{I}^2}^\8 \abs{\partial_x T_t(x,y)} &\frac{dt}{\sqrt{t}} dx  =  \int_{0}^{2^{-n}} \int_{2^{-2n} }^{xy \vee 2^{-2n}} ... + \int_{0}^{2^{-n}} \int_{xy \vee 2^{-2n}}^1 ...+ \int_{2^{n+2}}^{\8} \int_{2^{-2n}}^1...  \\
& + \int_{(2^{-n},2^{n+2}) \cap \set{|x-y|<2^{-n}}} \int_{2^{-2n} }^{1} ... +  \int_{(2^{-n},2^{n+2}) \cap \set{|x-y|>2^{-n}}} \int_{2^{-2n} }^{1} ...   \\
&+    \int_{0}^{\8} \int_{1}^{1 \vee \ln(\sqrt{xy})} ... +    \int_{0}^{\8} \int_{1 \vee \ln(\sqrt{xy})}^\8 ... \\
 =& A_{12} + A_{13}  + A_{14}  + A_{15} + A_{16}  + A_{17} + A_{18}.
}

For $A_{12}$ we have $ xy \sgeq t$, $t\leq 1$ and $x<y$, so that  $|F_2(t,x,y)| \sleq  y/t $. Using \eqref{R2}, \eqref{F2}, and \eqref{Th1},
\spx{
A_{12} & \sleq y \int_{0}^{2^{-n}} \int_{0}^{\8} \exp\left( -\frac{y^2}{ct} \right) \dtttt \, dx \\
& \sleq y^{-1}  \int_{0}^{2^{-n}} \, dx \cdot  \int_{0}^{\8} t^{-1} \exp\left( -\frac{1}{ct} \right) \, \dttt 
 \sleq 2^{-2n} \sleq 1.}

For $A_{13}$ we have $ xy \sleq t$, $t\leq 1$, and $x/t\sleq x^{-1}$, so that  $|F_1(t,x,y)| \sleq x^{-1} (xy/t)^{\be-1/2}$. Using \eqref{R1} and \eqref{F1},
\spx{A_{13} &  \sleq y^{\be} \cdot \int_{0}^{2^{-n}} x^{\be-1} \int_0^\8 t^{-\be} \exp\left( -\frac{y^2}{ct} \right) \frac{dt}{t} \, dx \\
& \sleq y^{-\be} \cdot \int_{0}^{2^{-n}} x^{\be-1}\, dx \cdot \int_0^\infty t^{-\be}  \exp\left( -\frac{1}{ct} \right) \frac{dt}{t}\\
& \sleq 2^{-2\be n}  \sleq 1,}
where in the last inequality we have used that $\be >0$.

For $A_{14}$ we have $ xy \sgeq t$, $|x-y|\simeq x$, and $x > y$, so that  $|F_2(t,x,y)| \sleq x/t$. Using \eqref{R2}, \eqref{Th1}, and \eqref{F2},
\spx{A_{14}  \sleq \int_{2^{n+2}}^{\infty} x \int_0^1 \exp\left( -\frac{x^2}{ct} \right) \dtttt \, dx \sleq \int_{2^{n+2}}^{\infty} x^{1-2N}  \, dx \cdot  \int_0^1 t^{N-2}\, dt \sleq 1.}


For $A_{15}$ we have that $ xy \sgeq t$, $x \simeq y \simeq 2^n$, and  $|F_3(t,x,y)| \sleq xt$. Using \eqref{R3}, \eqref{Th1}, and \eqref{F3},
\spx{A_{15}   \sleq &  \int_{\set{|x-y|<2^{-n}}} \int_{2^{-2n}}^1 \exp\eee{-cty^2} \left(\frac{|x-y|}{t} + xt \right) \frac{dt}{t} \, dx \\
 \sleq &  y^{-2N} \cdot  \int_{\set{|x-y|<2^{-n}}} |x-y| \, dx \cdot  \int_{2^{-2n}}^\8 t^{-N-2}\, dt \\
&+  y^{-2N} \cdot \int_{\set{|x-y|<2^{-n}}} x \, dx \cdot  \int_{2^{-2n}}^\8 t^{-N}\, dt \sleq 1.   }

For  $A_{16}$ we have that $xy \sgeq t$, $t\leq 1$, $x\sleq y$, and $\abs{F_3(t,x,y)} \sleq x^{-1} + xt$. Using \eqref{R3}, \eqref{Th1}, and \eqref{F3},
\spx{A_{16}   \sleq &  \int_{(2^{-n}, 2^{n+2}) \cap\set{|x-y|>2^{-n}}} \int_{2^{-2n}}^1 e^{-ty^2} \exp\eee{-\frac{|x-y|^2}{ct}} \eee{\frac{|x-y|}{t}+x^{-1} + tx} \frac{dt}{t} \, dx\\
= &A_{16,1}+ A_{16,2}+A_{16,3},   }
where  $A_{16,1}, A_{16,2}, A_{16,3}$ are the integrals with: $|x-y|t^{-1}, x^{-1}, xt$, respectively.
\spx{A_{16,1}  & \sleq  y^{-2N}  \int_{\set{|x-y|>2^{-n}}} |x-y| \int_{0}^\8 t^{-N-1}  \exp\eee{-\frac{|x-y|^2}{ct}} \, \frac{dt}{t} \, dx\\
&\sleq 2^{-2nN} \int_{\set{|x-y|>2^{-n}}} |x-y|^{-2N-1}\, dx \cdot  \int_{0}^\8 t^{-N-2}  \exp\eee{-\frac{1}{ct}} \, dt \sleq 1.
}
Notice that $x^{-1} \leq 2^n$, thus
\spx{A_{16,2}  & \sleq  y^{-2N} \int_{(2^{-n},\8) \cap \set{|x-y|>2^{-n}}} x^{-1} \int_{0}^\8 t^{-N}  \exp\eee{-\frac{|x-y|^2}{ct}} \, \frac{dt}{t} \, dx\\
&\sleq 2^{n(1-2N)} \int_{\set{|x-y|>2^{-n}}} |x-y|^{-2N}\, dx \cdot  \int_{0}^\8 t^{-N}  \exp\eee{-\frac{1}{ct}} \, \frac{dt}{t} \sleq 1.
}

\spx{A_{16,3}  & \sleq   \int_{0}^{2^{n+2}} x \, dx \cdot \int_{0}^\8 e^{-cty^2} \, dt  \sleq 2^{2n} \cdot 2^{-2n} \sleq 1.
}


For $A_{17}$ we have that $xy \sgeq \sh(2t)$, $t\geq 1$, and $\abs{F_2(t,x,y)} \sleq x+y \sleq y(x+1)$. Using \eqref{R2}, \eqref{Th2}, and \eqref{F2},
\spx{
A_{17}  \sleq & y e^{-cy^2} \int_0^\8 (x+1) e^{-cx^2} \, dx \cdot \int_1^\8 (\sh(2t))^{-1/2} \, \dtt \sleq 1.
}

For $A_{18}$ we have that $xy \sleq \sh(2t)$, $t\geq 1$, and $\abs{F_1(t,x,y)} \sleq  (xy/\sh(2t))^{\be-1/2} \cdot (x+x^{-1})$. Using \eqref{R1} and \eqref{F1},
\spx{
A_{18}   \sleq &  \int_0^{\8} \int_1^\8 e^{-c(x^2+y^2)} \eee{\frac{xy}{\sh(2t)}}^{\be} (x+x^{-1}) \, \frac{dt}{\sqrt{t\cdot\sh(2t)}} \, dx\\
\sleq & y^{\be} e^{-cy^2} \cdot \int_0^\8 x^{\be} (x+x^{-1}) e^{-cx^2} \, dx \cdot \int_1^\8 (\sh(2t))^{-\be-1/2} \, \dtt \sleq 1.
}


{\bf Proof of \eqref{a5} in Case 2.}
In this case we have $x,y \simeq 2^n$, $|x-y|\sleq 2^{-n} = d_I$. The proof follows by similar argument to those in  {\bf Case 1}. In particular, one uses \eqref{k123} and estimate $K_t^{[1]}$--$K_t^{[3]}$ in a similar way. The details are left to the reader.

{\bf Proof of \eqref{a6}.}  Let us write
\spx{
\int_X  \eee{x+x^{-1}} \int_0^\8 T_t(x,y) \dtt \, dx  = & \int_0^\8 \int_0^{1\land xy} ... + \int_0^\8 \int_{1\land xy}^1 ... \\
& + \int_{0}^\8 \int_1^{1\lor \ln(\sqrt{xy})} ... + \int_{0}^\8 \int_{1\lor \ln(\sqrt{xy})}^\8 ...\\
= & A_{19} + A_{20} + A_{21} + A_{22}.
}
Our goal is to prove $A_{19}+A_{20}+A_{21}+A_{22}\sleq 1$. Observe that by using \eqref{lag2-kernel}, \eqref{Th1}, \eqref{U_prop}, for $x\leq y/2$, we have
\sp{\label{rrr1}
\int_0^{1\land xy} T_t(x,y) \, \dtt \sleq \int_0^{xy} \exp\eee{-\frac{y^2}{ct}}  \dttt
\sleq 
\int_{y/x}^\8 e^{-ct} \, \dttt
\sleq e^{-cy/x}.
}
Similarly, we get the estimates
\al{
\label{rrr11}
&\int_0^{1\land xy} T_t(x,y) \, \dtt
\sleq e^{-cy^2}, && 2x\leq y, \ y\geq 1,\\
\label{rrr2}
&\int_0^{1\land xy} T_t(x,y) \, \dtt \sleq e^{-cx/y},  &&2x/3\geq y,\\
\label{rrr22}
&\int_0^{1\land xy} T_t(x,y) \, \dtt \sleq e^{-cx^2},  && 2x/3\geq y\geq 1.
}
Moreover, by \eqref{lag2-kernel}, \eqref{Th1}, \eqref{U_prop}, for $|x-y|\leq y/2$, we have
\sp{\label{rrr3}
\int_0^{1 \land xy} T_t(x,y) \, \dtt &\sleq \int_0^{ xy} \exp\eee{-\frac{|x-y|^2}{ct}}  \dttt\\
&\sleq 
\int_{|x-y|^2/(xy)}^\8 e^{-ct} \, \dttt \sleq 
 \ln\eee{\frac{y}{|x-y|}},
}
and, for $|x-y|\leq y/2$ and $y\geq 1$,
\sp{\label{rrr4}
\int_0^{1\land xy} T_t(x,y)  \, \dtt &\sleq \int_0^{1} \exp\eee{-\frac{|x-y|^2}{ct}}  \Theta(t,x,y) \dttt\\
&\sleq \int_0^{y^{-2}} \exp\eee{-\frac{|x-y|^2}{ct}} \dttt + \int_{y^{-2}}^1 (ty^2)^{-1} \exp\eee{-\frac{|x-y|^2}{ct}} \dttt \\
&\sleq \int_{y^{2}|x-y|^2}^\8 e^{-ct} \dttt + |x-y|^{-2}y^{-2} \int_{0}^{y^2|x-y|^2} e^{-ct} dt  \\
& \sleq \frac{\ln(2+(y|x-y|)^{-1})}{1+y^2|x-y|^2}.
}
Consider first $A_{19}$ in {\bf the case $y\leq 1$}. Using \eqref{rrr1}, \eqref{rrr3}, and \eqref{rrr2},
\spx{A_{19} \sleq & \int_0^{y/2}  e^{-cy/x}\, \frac{dx}{x} + y^{-1}\int_{|x-y|\leq y/2 } \ln\eee{\frac{y}{|x-y|}} \, dx \\
& + \int_{3y/2}^\8 xe^{-cx/y} \, dx +  \int_{3y/2}^\8 x^{-1} e^{-cx/y} \, dx
\sleq 1. }

Now consider $A_{19}$ in {\bf the case $y\geq 1$}. Using \eqref{rrr1}, \eqref{rrr11}, \eqref{rrr4}, and \eqref{rrr22}
\spx{A_{19} \sleq & \int_0^{(2y)^{-1}} e^{-cy/x} \, \frac{dx}{x} + \int_{(2y)^{-1}}^{y/2} (x+x^{-1})e^{-cy^2}\, dx \\&+ y\int_{|x-y|\leq y/2} \frac{\ln\eee{2+(y|x-y|)^{-1}}}{1+y^2|x-y|^2} \, dx
+ \int_{3y/2}^\8 x e^{-cx^2} \, dx \\
&\sleq e^{-cy}+ (y^2 + \ln y)e^{-cy^2} + \int_{-\8}^\8 \frac{\ln(2+|x|^{-1})}{1+x^2} dx+ \int_1^\8 x e^{-cx^2}\, dx \sleq 1. }

Recall that $\be > 0$. For $A_{20}$  we use \eqref{lag1-kernel} and \eqref{Bessel-function-small} getting
\spx{
A_{20} \sleq & \int_0^\8 (x+x^{-1}) \int_{0}^1 \eee{\frac{xy}{t}}^{\be} \exp\eee{-\frac{x^2+y^2}{ct}} \, \dttt\, dx\\
\sleq & \int_0^\8 (x+x^{-1}) \eee{\frac{xy}{x^2+y^2}}^{\be} \int_{x^2+y^2}^{\8} t^{\be}\exp(-t/c) \, \dttt \, dx\\
\sleq & \int_0^\8 (x+x^{-1}) \eee{\frac{xy}{x^2+y^2}}^{\be} \exp\eee{-c x^2} \, dx\\
 \sleq &\int_0^\8  \eee{\frac{xy}{x^2+y^2}}^{\be} \, \frac{dx}{x} \sleq 1.
}

For $A_{21}$ we have $xy \sgeq \sh(2t)$ and $x^{-1}\sleq y$ (otherwise $A_{21}=0$). Applying \eqref{lag2-kernel}, \eqref{Th2}, \eqref{U_prop}, we get
\spx{
A_{21} &\sleq \int_{0}^\8 (x+y) \int_1^{1\lor \ln(\sqrt{xy})} \sh(2t)^{-1/2} \Theta(t,x,y) \, \dtt \, dx\\
&\sleq (y+1) e^{-cy^2} \cdot \int_0^\8 (x+1) \exp(-cx^2) \, dx \cdot \int_1^\8 \sh(2t)^{-1/2} \, \dtt  \sleq 1.}

For $A_{22}$ we have $xy \sleq \sh(2t)$ and $\sh(2t)\simeq \ch(2t)$. Using \eqref{lag1-kernel} and \eqref{Bessel-function-small},
\spx{A_{22} & \sleq \int_{0}^\8 (x+x^{-1}) \int_{1\lor \ln(\sqrt{xy})}^\8 \sh(2t)^{-1/2} \eee{\frac{xy}{\sh(2t)}}^{\be} \exp\eee{-c(x^2+y^2)} \, \dtt \, dx \\
& \sleq y^{\be} e^{-cy^2} \cdot  \int_{0}^\8 (x+x^{-1}) x^{\be} e^{-cx^2} \, dx \cdot \int_1^\8 \sh(2t)^{-\be-1/2} \dtt \sleq 1.}
We have shown that $A_{19}+A_{20}+A_{21}+A_{22}\sleq 1$. This finishes the proofs of \eqref{a6} and Proposition \ref{a_lag_riesz}.

\subsection{Dirichlet Laplacian on $\RR^d_+$.} \label{sec:ex-dir-lap}

The goal of this section is to prove Theorem \ref{coro_dir}. Recall that we consider $X = \RR^{d}_+ = \RR^{d-1} \times (0,\8)$ and the Laplacian $L_D=-\Delta$ with the Dirichlet boundary condition at $x_d=0$. The semigroup generated by $-L_D$ is associated with the kernel
$$T_t(x,y) = H_t(\wt{x}, \wt{y}) T_{t,D}(x_d,y_d),$$
where $\wt{x} = (x_1,...,x_{d-1})$, $H_t$ is as in \eqref{heat_kernel}, and
\sp{\label{heat_kernel_dir}
T_{t,D}(x,y) = H_t(x-y) - H_t(x+y), \qquad x,y,t>0.
}
Observe that  $L_D$ is of the form as  in Section \ref{sec3.1}. Therefore, due to Theorem \ref{thm_riesz_loc_nonloc}, it is enough to prove Theorem \ref{coro_dir} in the case $d=1$.

Assume then that $d=1$, $X=(0,\8)$, and $L_D$ that is related to the kernel \eqref{heat_kernel_dir}. It is sufficient to check that for $T_{t,D}$ and the admissible covering
\sp{\label{dyadic}
\qq_D = \set{[2^n,2^{n+1}] \ : \ n \in \ZZ}
}
of $(0,\8)$  the conditions \eqref{a0}--\eqref{a5} are satisfied (obviously, in this case  \eqref{a6} is automatically satisfied). This is stated in the following lemma.

\lem{lem_dir_assumpt}{
Let $T_{t,D}(x,y)$ and $\qq_D$ be as in \eqref{heat_kernel_dir} and \eqref{dyadic}, respectively. Then \eqref{a0}--\eqref{a5} are satisfied.
}

\pr{[Sketch of the proof]
Notice that 
\spx{ \label{dirichlet-kernel-2}
T_{t,D}(x,y) = H_t(x,y) \Om(t,x,y), \qquad x,y,t>0,
}
where
\eqx{
\label{omega} 
\Om(t,x,y) = 1-\exp(-xy/t) \simeq \min(1,xy/t)
.
}

Moreover,
\sp{\label{der-dir}
\partial_x T_{t,D}(x,y) &= H_t(x,y) \eee{\Om(t,x,y)\frac{x-y}{t} + \frac{y}{t}\exp\left( -\frac{xy}{t}\right)} .
}

The conditions \eqref{a0}--\eqref{a5} can be proved by using standard estimates. For the convenience of the reader, we shall present the proof of \eqref{a4}. Let $I\in \qq_D$ and assume that $I=[2^n,2^{n+1}]$ for some $n\in\ZZ$.  Fix $y\in I^{**}$. From \eqref{der-dir} we have

We write
\spx{
\int_X \int_{d_I^2}^\8 \abs{\partial_x T_{t,D}(x,y)} \, \dtt \, dx & = \int_0^{2^{n+2}} \int_{d_I^2}^\8 ... +  \int_{2^{n+2}}^\8 \int_{d_I^2}^{\8} ...  \\
= A_1 + A_2.
}
For $A_1$ we have that $xy\sleq t$ and $x\sleq \sqrt{t}$. From \eqref{der-dir} we have $\abs{\partial_x T_{t,D}(x,y)} \sleq yt^{-3/2}(x/\sqrt{t}+1) \sleq   yt^{-3/2} $ and 
\eqx{
A_1\sleq \int_0^{2^{n+2}} \int_{2^{2n}}^\8 yt^{-3/2} \, \dtt \, dx \sleq 2^n \int_0^{2^{n+2}} \, dx \cdot  \int_{2^{2n}}^\8 t^{-2} \, dt \sleq 1.}
For $A_2$ we have that $|x|\simeq |x-y|$ and $y\sleq \sqrt{t} $. From \eqref{der-dir},
\spx{\label{der-dir2}
\abs{\partial_x T_{t,D}(x,y)} &\sleq t^{-1}  \exp\eee{-\frac{|x-y|^2}{ct}}\eee{\min(1,xy/t) + y/\sqrt{t}}\\
&\sleq t^{-1} \exp\eee{-\frac{x^2}{ct}} \cdot \begin{cases} \frac{y}{\sqrt{t}} & \text{ \ \ for } \frac{xy}{t} \sleq 1 \\ 1 & \text{ \ \  for }  \frac{xy}{t} \sgeq 1 \end{cases}.
}
Therefore,
\spx{
\int_{d_I^2}^\8 \abs{\partial_x T_{t,D}(x,y)}\dtt &\sleq \int_{d_I^2}^{xy} ... + \int_{xy}^{\8}... \\
&\sleq  \int_0^{xy} \exp\eee{-\frac{x^2}{ct}}\, \frac{dt}{t\sqrt{t}}  + y  \int_{0}^{\8} \exp\eee{-\frac{x^2}{ct}}\, \frac{dt}{t^2}  \\
&\sleq \int_0^{xy} \eee{\frac{t}{x^2}}^{3/2}\, \frac{dt}{t\sqrt{t}} +   yx^{-2} \sleq yx^{-2}.
}
Hence,
\eqx{A_2 \sleq 2^n \int_{2^{n+2}}^\8 x^{-2}\, dx \sleq 1.
}
This finishes the proof of \eqref{a4}.
}
\bibliographystyle{amsplain}        

\end{document}